\documentclass[12pt,reqno]{amsart}

\title[Optimal transport plans]
{Characterization of the optimal plans for the Monge-Kantorovich
 transport problem}
\author{Christian L\'eonard}
\date{July 2006}

\usepackage{amssymb, amsmath, amsfonts, latexsym, enumerate}

\newtheorem{theorem}[equation]{Theorem}
\newtheorem{lemma}[equation]{Lemma}
\newtheorem{proposition}[equation]{Proposition}
\newtheorem{corollary}[equation]{Corollary}

\newtheorem{definition}[equation]{Definition}

\theoremstyle{remark}
\newtheorem{remark}[equation]{Remark}
\newtheorem{remarks}[equation]{Remarks}

\newtheorem{examples}[equation]{Examples}

\numberwithin{equation}{section}

\newenvironment{listalpha}{%

  \begin{enumerate}}{\end{enumerate}}

\begin{document}


\newcommand{\R}{\mathbb{R}}

\newcommand{\1}{\textbf{1}}

\newcommand{\eqdef}{\stackrel{\vartriangle}{=}}
\newcommand{\dom}{\mathrm{dom\,}}
\newcommand{\icordom}{\mathrm{icordom\,}}
\newcommand{\cl}{\mathrm{cl\,}}
\newcommand{\cv}{\mathrm{cv}}
\newcommand{\inter}{\mathrm{int\,}}
\newcommand{\epi}{\mathrm{epi}\,}
\newcommand{\lsc}{lower semicontinuous}
\newcommand{\usc}{upper semicontinuous}
\newcommand{\cav}{{\hat{*}}}
\newcommand{\supp}{\mathrm{supp}\,}
\newcommand{\ls}{\mathrm{ls}\,}
\newcommand{\us}{\mathrm{us}\,}

\newcommand{\boulette}[1]{$\bullet$\ Proof of (#1).}
\newcommand{\Boulette}[1]{\par\medskip\noindent $\bullet$\ Proof of (#1).}

\newcommand{\limn}{\lim_{n\rightarrow\infty}}

\newcommand{\positif}{positif}
\newcommand{\negatif}{n\'egatif}
\newcommand{\croissant}{croissant}
\newcommand{\decroissant}{d\'ecroissant}
\newcommand{\strictpositif}{strictly positif}
\newcommand{\strictnegatif}{strictly n\'egatif}
\newcommand{\strictcroissant}{strictly croissant}
\newcommand{\strictdecroissant}{strictly d\'ecroissant}

\newcommand{\aldu}{\ast}
\newcommand{\diffdom}{\mathrm{diffdom\,}}

\newcommand{\UU}{\mathcal{U}}
\newcommand{\LL}{\mathcal{L}}
\newcommand{\YY}{\mathcal{Y}}
\newcommand{\XX}{\mathcal{X}}
\newcommand{\Li}{\mathcal{L}_1}
\renewcommand{\Xi}{\mathcal{X}_1}
\newcommand{\Ui}{\mathcal{U}_1}
\newcommand{\Uii}{\Li^\aldu}
\newcommand{\Yi}{\mathcal{Y}_1}
\newcommand{\Yii}{\Xi^\aldu}
\newcommand{\Ci}{C_1}
\newcommand{\Ss}{S^{\aldu}}
\newcommand{\Ct}{C_\theta}
\newcommand{\Cts}{C_{\theta^*}}
\newcommand{\AB}{{A\!\times\! B}}
\newcommand{\PA}{\mathcal{P}_A}
\newcommand{\PB}{\mathcal{P}_B}
\newcommand{\PAB}{\mathcal{P}_{AB}}
\newcommand{\CAB}{C_{AB}}
\newcommand{\CAs}{C_A^*}
\newcommand{\CBs}{C_B^*}
\newcommand{\CABs}{\CAB^*}
\newcommand{\li}{B_b}
\newcommand{\lip}{B_b.\pi}
\newcommand{\Pc}{\mathcal{P}_c}
\newcommand{\Q}{\mathcal{Q}}
\newcommand{\Qp}{\Q^+}
\newcommand{\Lp}{\mathcal{L}^+}
\newcommand{\Mp}{\mathcal{M}^+_c}
\newcommand{\M}{\mathcal{M}_c}
\newcommand{\s}{\mathcal{S}}
\newcommand{\sA}{{\mathcal{S}_A}}
\newcommand{\sB}{{\mathcal{S}_B}}
\newcommand{\PS}{\mathcal{P}_{\s}}
\newcommand{\Cc}{C_{c}}

\newcommand{\F}{\Phi}
\newcommand{\Fs}{\Phi^*}
\newcommand{\Fss}{\Phi^{**}}
\newcommand{\Fd}{\Phi_2}
\newcommand{\Fb}{\overline{\Phi}}
\newcommand{\La}{\Lambda}
\newcommand{\Ls}{\Lambda^*}
\newcommand{\Lb}{\Lambda_2}
\newcommand{\fss}{\varphi^{**}}
\newcommand{\ts}{\theta^*}
\newcommand{\jt}{j_{\theta}}
\newcommand{\jts}{j_{\theta^*}}
\newcommand{\ds}{\delta_{\Ct}^*}

\newcommand{\Po}{$(P)$}
\renewcommand{\Pi}{$(P_1)$}
\newcommand{\PX}{$(P_{1,\XX})$}
\newcommand{\Pxb}{$(P_1^{\xb})$}
\newcommand{\Do}{$(D_0)$}
\newcommand{\Di}{$(D_1)$}
\newcommand{\Dii}{$(D_2)$}
\newcommand{\Dxb}{$(D_2^{\xb})$}

\newcommand{\NF}{|\cdot|_\F}
\newcommand{\NL}{|\cdot|_\La}
\newcommand{\sLU}{\sigma(\LL,\UU)}
\newcommand{\sLUi}{\sigma(\Li,\Ui)}
\newcommand{\sUL}{\sigma(\UU,\LL)}
\newcommand{\sULi}{\sigma(\Ui,\Li)}
\newcommand{\sULii}{\sigma(\Uii,\Li)}
\newcommand{\sXY}{\sigma(\XX,\YY)}
\newcommand{\sXYi}{\sigma(\Xi,\Yi)}
\newcommand{\sXYii}{\sigma(\Xi,\Yii)}
\newcommand{\sYXii}{\sigma(\Yii,\Xi)}
\newcommand{\sYX}{\sigma(\YY,\XX)}
\newcommand{\sYXi}{\sigma(\Yi,\Xi)}
\newcommand{\sYXx}{\sYXii}

\newcommand{\HF}{$(H_\F)$}
\newcommand{\HFi}{$(H_{\F1})$}
\newcommand{\HFii}{$(H_{\F2})$}
\newcommand{\HFiii}{$(H_{\F3})$}
\newcommand{\HT}{$(H_T)$}
\newcommand{\HTi}{$(H_{T1})$}
\newcommand{\HTii}{$(H_{T2})$}
\newcommand{\HC}{$(H_{C})$}

\newcommand{\ul}{\langle u,\ell\rangle}
\newcommand{\lz}{\langle \ell,\zeta\rangle}
\newcommand{\xo}{\langle x,\omega\rangle}
\newcommand{\yx}{\langle y,x\rangle}
\newcommand{\xy}{\langle XY\rangle}
\newcommand{\lb}{\bar\ell}
\newcommand{\ub}{\bar u}
\newcommand{\ob}{\bar \omega}
\newcommand{\yb}{\bar y}
\newcommand{\xb}{\bar x}
\newcommand{\bq}{\langle b,q\rangle}
\newcommand{\ap}{\langle a,p\rangle}
\newcommand{\ab}{\bar a}
\newcommand{\bb}{\bar b}
\newcommand{\IAB}{\int_{\AB}}
\newcommand{\IS}{\int_{\s}}
\newcommand{\uS}{u_{|\s}}
\newcommand{\etac}{\check{\eta}}
\newcommand{\etat}{\widetilde{\eta}}
\newcommand{\thetat}{\widetilde{\theta}}
\newcommand{\etab}{\overline{\eta}}
\newcommand{\thetab}{\overline{\theta}}
\newcommand{\etah}{\widehat{\eta}}
\newcommand{\thetah}{\widehat{\theta}}
\newcommand{\suppp}{\mathrm{supp\,}\pi}

\newcommand{\IZ}{\int_{Z}}
\newcommand{\Co}{C}
\newcommand{\Uo}{\UU}
\newcommand{\Lo}{\LL}
\newcommand{\lmax}{{\lambda_o}}
\newcommand{\lmaxs}{{\lambda_o^*}}
\newcommand{\Qb}{\overline{Q}}
\newcommand{\MZ}{\mathrm{Meas}(Z)}
\newcommand{\MU}{\mathrm{Meas}_{\lmax}(Z)}
\newcommand{\Lls}{L_{\lmaxs}}
\newcommand{\MK}{M\! K}


 \address{Modal-X, Universit\'e Paris 10\quad \& }
 \address{CMAP, \'Ecole Polytechnique. 91128 Palaiseau Cedex, France}
 \email{christian.leonard@polytechnique.fr}
 \keywords{Convex optimization, saddle-point, conjugate duality, optimal transport}
 \subjclass[2000]{46N10, 49J45, 28A35}

\begin{abstract}
We present a general method, based on conjugate duality, for
solving a convex minimization problem without assuming unnecessary
topological restrictions on the constraint set. It leads to dual
equalities and characterizations of the minimizers without
constraint qualification.
\\
As an example of application, the Monge-Kantorovich optimal
transport problem is solved in great detail. In particular, the
optimal transport plans are characterized without restriction.
This characterization improves the already existing literature on
the subject.
\end{abstract}

\maketitle
\tableofcontents


\section{Introduction}\label{sec:introduction}

Although the title highlights Monge-Kantorovich optimal transport
problem, the aim of this paper is twofold.
\begin{itemize}
    \item First, one presents an \emph{``extended'' saddle-point method} for
solving a convex minimization problem: It is shown how to
implement the standard saddle-point method in such a way that
topological restrictions on the constraint sets (the so-called
constraint qualifications) may essentially be removed. Of course,
so doing one has to pay the price of solving an arising new
problem. Namely, one has to compute the extension of some
function; this may be a rather difficult task in some situations,
but it will be immediate in the Monge-Kantorovich case. This
method is based on conjugate duality as developed by
R.T.~Rockafellar in \cite{Roc74}. Dual equalities and
characterizations of the minimizers are obtained without
constraint qualification.
    \item Then, these ``extended'' saddle-point abstract results are applied
to the Monge-Kantorovich optimal transport problem. In particular,
the optimal plans are characterized without any restriction. This
characterization improves the already existing literature on the
subject.
\end{itemize}

Other applications of the extended saddle-point method are
investigated by the author in \cite{Leo06b} in connection with
entropy minimization.

\subsection*{The Monge-Kantorovich transport problem}

Let us take $A$ and $B$ two Polish  (separable complete metric)
spaces furnished with their respective Borel $\sigma$-fields, a
\lsc\ (cost) function $c:\AB \to [0,\infty]$ which may take
infinite values and two probability measures $\mu\in\PA$ and
$\nu\in\PB$ on $A$ and $B.$ We denote $\PA, \PB$ and $\PAB$ the
sets of all Borel probability measures on $A,$ $B$ and $\AB.$ The
Monge-Kantorovich problem is
\begin{equation}\label{MK}
    \textsl{minimize } \pi\in\PAB\mapsto \IAB c(a,b)\,\pi(dadb)
    \textsl{ subject to }\pi\in P(\mu,\nu) \tag{$\MK$}
\end{equation}
where $P(\mu,\nu)$ is the set of all $\pi\in\PAB$ with prescribed
marginals $\pi_A=\mu$ on $A$ and $\pi_B=\nu$ on $B.$ Note that $c$
is measurable since it is \lsc\ and the integral $\IAB
c\,d\pi\in[0,\infty]$ is well-defined since $c\geq 0.$
\\
For a general account on this active field of research, see the
books of S.~Rachev and L.~R\"uschendorf \cite{RacRus} and
C.~Villani \cite{Vill03,Vill05}.

\begin{definition}[Optimal plan]\label{def-02}
One says that  $\pi\in P(\mu,\nu)$ is an \emph{optimal plan} if it
minimizes $\gamma\mapsto\IAB c\,d\gamma$ on $P(\mu,\nu)$
\emph{and} $\IAB c\,d\pi<\infty.$
\end{definition}

It is well-known that there exists at least an optimal plan if and
only if there exists some $\pi^o\in P(\mu,\nu)$  such that $\IAB
c\,d\pi^o<\infty;$ this will be recovered at Theorem
\ref{res-MK1}. Definition \ref{def-02} throws away the
uninteresting case where $\IAB c\,d\pi=\infty$ for all $\pi\in
P(\mu,\nu).$ Note also that, since Monge-Kantorovich problem is
not a strictly convex problem, infinitely many optimal plans may
exist.

\subsubsection*{Already existing optimality criteria}

Some usual criteria are expressed in terms of cyclical
$c$-monotonicity.
\begin{definition}[Cyclically $c$-monotone plan]
A subset $\Gamma\subset\AB$ is said to be
\emph{cyclically $c$-monotone} if for any integer $n\geq 1$ and
any family $(a_1,b_1),\dots,(a_n,b_n)$ of points in $\Gamma,$
 $
    \sum_{i=1}^n c(a_i,b_i)\leq \sum_{i=1}^n c(a_i,b_{i+1})
 $
with the convention $b_{n+1}=b_1.$
\\
A probability measure $\pi\in\PAB$ is said to be cyclically
$c$-monotone if it is concentrated on a measurable cyclically
$c$-monotone set $\Gamma,$ i.e. $\pi(\Gamma)=1.$
\end{definition}
This notion goes back to the seminal paper \cite{Rus96} by
L.~R\"uschendorf where the standard cyclical monotonicity of
convex functions introduced by Rockafellar has been extended in
view of solving Monge-Kantorovich problem.
\\
While completing this paper, the author has been informed of the
recent work \cite{ST06} by W.~Schachermayer and J.~Teichman who
have improved previous characterization criteria  in several
directions.  The following definition introduced in \cite{ST06} is
useful to state \cite{ST06}'s results in a concise way.
\begin{definition}[Strongly $c$-monotone plan]
A transport plan $\pi\in P(\mu,\nu)$ is called \emph{strongly
$c$-monotone} if there exist two measurable functions $\varphi$
and
    $\psi$ on $A$ and $B$ taking their values in
    $[-\infty,+\infty)$ such that
    \begin{equation}\label{eq-75}
    \left\{\begin{array}{l}
      \varphi\oplus\psi\leq c\quad \textrm{everywhere} \\
       \varphi\oplus\psi= c\quad \pi\textrm{-almost everywhere.} \\
    \end{array}\right.
    \end{equation}
\end{definition}
Here and below, we denote
$\varphi\oplus\psi(a,b)=\varphi(a)+\psi(b).$
\\
One easily shows that a strongly $c$-monotone plan is cyclically
$c$-monotone.
\\
The main results of \cite{ST06} are collected in the next two
theorems.
\begin{theorem}[\cite{ST06}]\label{res-15}
Let $c$ be a \lsc\ nonnegative \emph{finitely-valued} function. If
there exists some $\pi^o\in P(\mu,\nu)$ such that $\IAB
c\,d\pi^o<\infty,$ then for any $\pi\in P(\mu,\nu),$ the following
three statements are equivalent:
\begin{enumerate}[(i)]
    \item $\pi$ is an optimal plan;
    \item $\pi$ is cyclically $c$-monotone;
    \item $\pi$ is strongly $c$-monotone.
\end{enumerate}
\end{theorem}
This result significantly improves an already existing criterion
(see \cite{Vill05}, Chapter 5) where the same conclusion holds
with a finitely-valued function $c$ under the following constraint
qualification: There exist two nonnegative measurable functions
$c_A$ and $c_B$ on $A$ and $B$ such that
\begin{equation}\label{eq-29a}
    c\leq c_A\oplus c_B, \int_A c_A\,d\mu<\infty \textrm{ and }
    \int_B c_B\,d\nu<\infty.
\end{equation}
Note that (\ref{eq-29a}) implies that $\IAB c\,d\pi<\infty$ for
\emph{all} $\pi\in P(\mu,\nu).$ It also improves a result of
L.~Ambrosio and A.~Pratelli \cite{AP02} who have shown that, when
$c$ is finitely-valued and under the moment condition
\begin{equation}\label{eq-29b}
\left.\begin{array}{l}
  \mu\left(\left\{a\in A; \int_B c(a,b)\,\nu(db)<\infty\right\}\right)>0 \\
  \nu\left(\left\{b\in B; \int_A c(a,b)\,\mu(da)<\infty\right\}\right)>0
\end{array}
\right.
\end{equation}
 which is weaker than (\ref{eq-29a}), any cyclically $c$-monotone $\pi$ in $P(\mu,\nu)$ is both an
 optimal and a strongly $c$-monotone plan.
 For (\ref{eq-29b}) to hold, it is enough that $\IAB
c\,d\mu\otimes\nu<\infty.$ It is also proved in \cite{AP02} that
the functions $\varphi$ and $\psi$ in (\ref{eq-75}) can be taken
such that $\varphi\in L_1(A,\mu)$ and $\psi\in L_1(B,\nu).$

The next result is concerned with cost functions $c$ which may
take infinite values.
\begin{theorem}[\cite{AP02,ST06}]
Let $c$ be a \lsc\ $[0,\infty]$-valued function.
\begin{enumerate}[(a)]
    \item Any optimal plan is  cyclically $c$-monotone.
    \item If
    \begin{equation}\label{eq-76}
    \mu\otimes\nu(\{c< \infty\})=1,
    \end{equation}
   then any optimal plan is strongly $c$-monotone.
    \item If there exists some $\pi^o\in P(\mu,\nu)$ such that $\IAB
    c\,d\pi^o<\infty,$ then any strongly $c$-monotone plan in $P(\mu,\nu)$
    is an optimal plan.
\end{enumerate}
\end{theorem}
Statement (a) is proved in \cite{AP02}, while statements (b) and
(c) are taken from \cite{ST06}.

\begin{examples}\
\begin{enumerate}
    \item An interesting example of a cyclically $c$-monotone plan which is
not optimal is exhibited in \cite{AP02}, in a situation where $c$
takes infinite values and an optimal plan exists. This is in
contrast with Theorem \ref{res-15} and emphasizes that cyclical
$c$-monotonicity isn't the right notion to consider in the general
case.
    \item Take $A=B=[0,1],$ $\mu(da)=da,$ $\nu(db)=db$ the Lebesgue measure
on $[0,1]$ and $c(a,b)=0$ if $a=b$ and $+\infty$ otherwise.
Condition (\ref{eq-76}) is restrictive enough to rule this basic
situation out. In the present paper, this restriction is removed.
\end{enumerate}
\end{examples}

\subsubsection*{A new optimality criterion}

Our main results about the optimal plans are Theorems
\ref{res-MK2} and \ref{res-MK3}. Next theorem sums them up.

\begin{theorem}\label{res-21}
Let $c$ be a \lsc\ $[0,\infty]$-valued function and let $\pi\in
P(\mu,\nu)$ satisfy $\IAB c\,d\pi<\infty.$
\begin{enumerate}[(a)]
    \item
    $\pi$ is an optimal plan if and only if there exist
     two finitely-valued functions $\varphi\in\mathbb{R}^A$ and $\psi\in\mathbb{R}^B$  such that
    \begin{equation}\label{eq-75b}
    \left\{\begin{array}{l}
   \varphi\oplus\psi\leq c \quad \textrm{everywhere and}\\
    \varphi\oplus\psi= c \quad \pi\textrm{-almost everywhere.} \\
    \end{array}\right.
    \end{equation}

    \item
If  $\pi$ is an optimal plan, there exist two finitely-valued
functions $\varphi\in\mathbb{R}^A$ and $\psi\in\mathbb{R}^B$  such
that $\varphi\in L_1(A,\mu),$ $\psi\in L_1(B,\nu),$
\begin{equation}\label{eq-75c}
    \left\{\begin{array}{clll}
      |\varphi\oplus\psi|&\leq &c & \textrm{everywhere and} \\
       \varphi\oplus\psi&= &c &  \textrm{on }\suppp\bigcap\{c<\infty\}.\\
    \end{array}\right.
\end{equation}
\end{enumerate}
\end{theorem}

\begin{remarks} These results improve previous literature on the
subject in several aspects.
\begin{enumerate}[a.]
\item No restriction is imposed on $c,$ $\mu$ and $\nu.$ In
particular, (\ref{eq-76}) is removed.

\item For the optimality criterion (a), the so-called Kantorovich
potentials $\varphi$ and $\psi$ are finitely-valued and are not
required to be a priori measurable. This is in contrast with the
definition of strongly $c$-monotone plans.

\item The analogue of (b) is usually stated as follows: If  $\pi$
is an optimal plan, there exist two $[-\infty,\infty)$-valued
functions $\varphi\in L_1(A,\mu)$ and $\psi\in L_1(B,\nu)$ such
that (\ref{eq-75b}) holds, even in the case where $c$ is required
to be finite. The improvements carried by (\ref{eq-75c}) are:
\begin{itemize}
    \item[-] The equality $ \varphi\oplus\psi=c $ holds on
    $\suppp\cap\{c<\infty\}$ rather than only $\pi$-almost
    everywhere;
    \item[-] The Kantorovich potentials $\varphi$ and $\psi$ are
    finitely-valued;
    \item[-] We obtain $|\varphi\oplus\psi|\leq c$ rather than $\varphi\oplus\psi\leq c.$
\end{itemize}
\end{enumerate}
\end{remarks}
As an immediate consequence of Theorem \ref{res-21}, we obtain the
following
\begin{corollary}
Any $\pi\in P(\mu,\nu)$ satisfying $\IAB c\,d\pi<\infty$ is an
optimal plan if and only if it is strongly $c$-monotone.
\end{corollary}
But the sufficient condition of Theorem \ref{res-21} is weaker
than the strong $c$-monotonicity, while its necessary condition is
stronger.

Finally, let us indicate why considering cost functions $c$
possibly achieving the value $+\infty$ is a significant extension.
In the finite-valued case, the domain of $c$ is the closed
rectangle $\AB.$ If one wants to forbid transporting mass from $A$
to $B$ outside some closed subset $\s$ of $\AB$ and only consider
the finitely-valued \lsc\ cost function $\tilde{c}$ on $\s,$
simply consider the extended cost function $c$ on $\AB$ which
matches with $\tilde{c}$ on $\s$ and is $+\infty$ outside. In this
case, $c$ has a closed effective domain. But there are also \lsc\
functions $c$ whose domain is an increasing union of closed
subsets.

\subsection*{An abstract convex problem and related questions}
Monge-Kantorovich problem is a particular instance of an abstract
convex minimization problem which we present now.

Let $\UU$ be a vector space, $\LL=\UU^\aldu$ its algebraic dual
space, $\F$ a $(-\infty,+\infty]$-valued convex function on $\UU$
and $\Fs$ its convex conjugate for the duality $\langle
\UU,\LL\rangle.$ Let $\YY$ be another vector space,
$\XX=\YY^\aldu$ its algebraic dual space and $T:\LL\rightarrow\XX$
is a linear operator. We consider the convex minimization problem
\begin{equation}
  \textsl{minimize } \Fs(\ell) \textsl{ subject to } T\ell\in C,\
  \ell\in\LL
\tag{$P$}
\end{equation}
where $ C$ is a convex subset of $\XX.$ As is well known,
Fenchel's duality leads to the dual problem
\begin{equation}
  \textsl{maximize } \inf_{x\in C}\yx - \F(T^\ast y), y\in\YY
\tag{$D$}
\end{equation}
where $T^\ast $ is the adjoint of $T.$

What about Monge-Kantorovich problem?  We denote $C_A,$ $C_B$ and
$\CAB$ the spaces of all continuous bounded functions on $A,$ $B$
and $\AB;$ $C_A^*,$ $C_B^*$ and $\CAB^*$ are their algebraic dual
spaces. Taking $\LL=\CAB^*$  the algebraic dual of $\UU=C_{AB},$
$T$ will be the marginal operator $T\ell=(\ell_A,\ell_B)\in
\XX:=C_A^*\times C_B^*$ which in restriction to those $\ell$'s in
$\LL$ which are probability measures gives the marginals $\ell_A$
on $A$ and $\ell_B$ on $B$ and $C$ will simply be $\{(\mu,\nu)\}.$
Choosing $\F(\varphi,\psi)=0$ if $\varphi\oplus \psi\leq c$ and
$\F(\varphi,\psi)=+\infty$ otherwise, will lead us to
Monge-Kantorovich problem.

The usual questions related to $(P)$ and $(D)$ are
\begin{itemize}
 \item the dual equality: Does $\inf(P)=\sup(D)$ hold?
 \item the primal attainment: Does there exist a solution $\lb$ to
    $(P)?$ What about the minimizing sequences, if any?
 \item the dual attainment: Does there exist a solution $\yb$ to $(D)?$
 \item the representation of the primal solutions: Find an identity
    of the type: $\lb\in\partial\F(T^\ast \yb).$
\end{itemize}
 We are going to answer them in terms of
some \emph{extension} $\bar\Phi$ of $\F$ under the weak assumption
\begin{equation}\label{eq-14b}
    T^{-1}(C)\cap\diffdom\Fs\not =\emptyset
\end{equation}
where
 $$
\diffdom\Fs=\{\ell\in\LL;\partial_{\LL^\aldu}\Fs(\ell)\not=\emptyset\}
 $$
is the subset of all vectors in  $\LL$ at which $\Fs$ admits a
nonempty subdifferential with respect to the algebraic dual
pairing $\langle\LL,\LL^\aldu\rangle$ where $\LL^\aldu$ is the
algebraic dual space of $\LL.$ Note that by the geometric version
of Hahn-Banach theorem, the intrinsic core of the effective domain
of the objective function $\Fs:$ $\icordom\Fs,$ is included in
$\diffdom\Fs.$ Hence, a useful criterion to get (\ref{eq-14b}) is
\begin{equation}\label{eq-14}
    T^{-1}(C)\cap\icordom\Fs\not =\emptyset.
\end{equation}
The drawback of such a general approach is that one has to compute
the extension $\bar\Phi.$ In specific examples, this might be a
difficult task. The extension $\bar\Phi$ is made precise at
Section \ref{sec:transport} for Monge-Kantorovich problem. Another
important example of application of our general results is the
problem of minimizing an entropy functional under a convex
constraint. This is worked out by the author in \cite{Leo06b} with
probabilistic applications in mind; it is based on the explicit
expression of the corresponding function $\bar\Phi.$

The restriction (\ref{eq-14}) seems very weak since $\icordom\Fs$
is the notion of interior which gives the largest possible set. As
$T^{-1}(C)\cap\dom\Fs=\emptyset$ implies that $(P)$ has no
solution, the only case where the problem remains open when
$\icordom\Fs$ is nonempty is the situation where $T^{-1}(C)$ and
$\dom\Fs$ are tangent to each other. This is used in \cite{Leo06b}
to obtain general results for convex integral functionals.
\\
Nevertheless, the Monge-Kantorovich optimal transport problem
provides an interesting case where the  constraints never stand in
$\icordom\Fs$ (see Remark \ref{rem-04}) so that (\ref{eq-14}) is
useless and (\ref{eq-14b}) is the right assumption to be used.

\subsubsection*{The strategy} A usual way to prove the dual
attainment and obtain some representation of the primal solutions
is to require that the constraint is qualified: a property which
allows to separate the convex constraint set $T^{-1}(C)$ and the
level sets of the objective function. The strategy of this article
is different: one chooses ad hoc topologies so that the level sets
have nonempty interiors. This also allows to apply Hahn-Banach
theorem, but this time the constraint set is not required to be
qualified. We take the rule not to introduce arbitrary topological
assumptions since $(P)$ is expressed without any topological
notion. Because of the convexity of the problem, one takes
advantage of geometric easy properties: the topologies to be
considered later are associated with seminorms which are gauges of
level sets of the convex functions $\F$ and $\Fs.$ They are useful
tools to work with the \emph{geometry} of $(P).$
\\
It appears that when the constraints are infinite-dimensional one
can choose several different spaces $\YY$ without modifying the
value and the solutions of $(P).$ So that for a small space $\YY$
the dual attainment is not the rule. As a consequence, we are
facing the problem of finding an \emph{extension} of $(D)$ which
admits solutions in generic cases and such that the representation
of the primal solution is $\lb\in\partial\bar\Phi(T^\ast \yb)$
where $\bar\Phi$ is some extension of $\F.$
\\
We are going to
\begin{itemize}
    \item use the standard saddle-point approach to convex problems
    based on conjugate duality as developed by Rockafellar in \cite{Roc74}
    \item with topologies which reflect some of the geometric
    structure of the objective function.
\end{itemize}
These made-to-measure topologies are associated with the gauges of
the level sets of $\F$ and $\Fs.$

\subsection*{Outline of the paper}
The abstract results are stated without proof at Section
\ref{sec:abstractpb}. Their proofs are postponed to Section
\ref{sec:proofs}. Section \ref{sec:transport} is devoted to the
application of the abstract results to the Monge-Kantorovich
optimal transport problem. Finally, basic results about convex
minimization and gauge functionals are recalled in the Appendix.

\subsection*{Notation}
Let $X$ and $Y$ be topological vector spaces. The algebraic dual
space of $X$ is $X^{\aldu},$ the topological dual space of $X$ is
$X'.$ The topology of $X$ weakened by $Y$ is $\sigma(X,Y)$ and one
writes $\langle X,Y\rangle$ to specify that $X$ and $Y$ are in
separating duality.
\\
Let $f: X\rightarrow [-\infty,+\infty]$ be an extended numerical
function. Its convex conjugate with respect to $\langle
X,Y\rangle$ is $f^*(y)=\sup_{x\in X}\{\langle x,y\rangle
-f(x)\}\in [-\infty,+\infty],$ $y\in Y.$ Its subdifferential at
$x$ with respect to $\langle X,Y\rangle$ is $\partial_Y
f(x)=\{y\in Y; f(x+\xi)\geq f(x)+\langle y,\xi\rangle, \forall
 \xi\in X\}.$ If no confusion occurs, one writes
$\partial f(x).$
\\
For each point $a,$ $\epsilon_a$ is the Dirac measure at $a.$

\subsection*{Stop saying no, be strict}\footnote{This is only a suggestion, not a demanding of the right to
write maths differently.} The function $\sin x$ is not negative,
but it is not nonnegative. It is not decreasing, but it is not
nondecreasing. All this does not make much sense and is not far
from being a nonsense for non-English speaking people. As a
convention, we'll use the non-English way of saying that a
\emph{positif} function is a $[0,\infty)$-valued function while if
it is $(0,\infty)$-valued it is also \emph{strictly positif}. The
integer part is a \emph{croissant} (increasing in colloquial
English) function and the exponential is also a \emph{strictly
croissant} function.  Symmetrically, we also use the notions of
\emph{n\'egatif} (negative in colloquial English) and strictly
n\'egatif, \emph{d\'ecroissant} (decreasing in colloquial English)
and strictly d\'ecroissant functions or sequences. To be coherent,
$[0,\infty)$ and $(-\infty,0]$ are respectively the sets of
positif and n\'egatif numbers, and $\epsilon >0$ is also strictly
positif. We keep the French words not to be mixed up with the
usual way of writing mathematics in English.

\section{The abstract convex minimization problem}
\label{sec:abstractpb}

In this section we give the statements of the results about the
abstract convex minimization problem. The dual equality and the
primal attainment are stated at Theorem \ref{xP3}; the dual
attainment and the dual representation of the minimizers are
stated at Theorems \ref{T3a} and \ref{T3b}. Their proofs are
postponed to Section \ref{sec:proofs}.

\subsection{Basic diagram}
\label{sec:Tadj}

Let $\YY$ be a vector space and $\XX=\YY^\aldu$ its algebraic dual
space. It is useful to define the constraint operator $T$ by means
of its adjoint $T^\ast :\YY\rightarrow\LL^\aldu$ $(\LL^\aldu$ is
the algebraic dual space of $\LL),$ as follows. For all
$\ell\in\LL, x\in\XX,$
\begin{displaymath}
  T\ell= x
  \Longleftrightarrow
  \forall y\in\YY, \langle
  T^\ast y,\ell\rangle_{\LL^\aldu,\LL}= \yx_{\YY,\XX}.
\end{displaymath}
We shall assume that the restriction
\begin{equation}\label{eq-01}
    T^\ast (\YY)\subset \UU
\end{equation}
 holds, where $\UU$ is identified with a subspace of
$\LL^\aldu=\UU^{\aldu\aldu}.$ It follows that the diagram
\begin{equation}
 \begin{array}{ccc}
\Big\langle\ \UU & , & \LL \ \Big\rangle \\
T^\ast  \Big\uparrow & & \Big\downarrow
 T
\\
\Big\langle\ \YY & , & \XX\ \Big\rangle
\end{array}
\tag{Diagram 0}
\end{equation}
is meaningful.

\subsection{Assumptions}
\label{sec:assumptions}

Let us give the list of our main hypotheses.
\begin{itemize}
\item[\HF]
1-\quad $\F: \UU\rightarrow [0,+\infty]$ is convex and $\F(0)=0$\\
2-\quad $\forall u\in\UU, \exists \alpha >0, \F(\alpha u)<\infty$\\
3-\quad $\forall u\in\UU, u\not=0,\exists t\in\R, \F(tu)>0$
\item[\HT]
1-\quad $T^\ast (\YY)\subset\UU$\\
2-\quad $\mathrm{ker\ }T^\ast =\{0\}$
\item [\HC] \qquad
$\Ci\eqdef
 C\cap\Xi$ is a convex $\sigma(\Xi,\Yi)$-closed subset of $\Xi$
\end{itemize}
The definitions of the vector spaces $\Xi$ and $\Yi$ which appear
in the last assumption are stated below at Section
\ref{sec:problems}. For the moment, let us only say that if $ C$
is convex and $\sigma(\XX,\YY)$-closed, then \HC\ holds.

\par\smallskip
\noindent\textit{Comments about the assumptions.}\
\begin{itemize}
    \item[-] By construction, $\Fs$ is a convex $\sigma(\LL,\UU)$-closed
function, even if $\F$ is not convex. Assuming the convexity of
$\F$ is not a restriction.
    \item[-] The assumption \HFi\ also expresses that $\F$ achieves its minimum
at $u=0$ and that $\F(0)=0.$ This is a practical normalization
requirement which will allow us to build a gauge functional
associated with $\F.$ More, \HFi\ implies that $\Fs$ also shares
this property. Gauge functionals related to $\Fs$ will also appear
later.
    \item[-] With any convex function $\tilde\F$ satisfying \HFii, one can
associate a function $\F$ satisfying \HFi\ in the following
manner. Because of \HFii, $\tilde\F(0)$ is finite and there exists
$\ell_o\in\LL$ such that $\ell_o\in\partial\tilde\F(0).$ Then,
$\F(u)\eqdef\tilde\F(u)-\langle\ell_o,u\rangle -\tilde\F(0),$
$u\in\UU,$ satisfies \HFi\ and
${\tilde\F}^*(\ell)=\Fs(\ell-\ell_o)-\tilde\F(0),$ $\ell\in\LL.$
    \item[-] The hypothesis \HFiii\ is not a restriction. Indeed, assuming \HFi, let us
suppose that there exists a direction $u_o\not=0$ such that
$\F(tu_o)=0$ for all real $t.$ Then any $\ell\in\LL$ such that
$\langle\ell,u_o\rangle\not=0$ satisfies
$\Fs(\ell)\geq\sup_{t\in\R}t\langle\ell,u_o\rangle=+\infty$ and
can't be a solution to $(P).$
    \item[-] The hypothesis \HTii\ isn't a restriction either: If
$y_1-y_2\in\mathrm{ker\ }T^\ast ,$ we have $\langle
T\ell,y_1\rangle=\langle T\ell,y_2\rangle,$ for all $\ell\in\LL.$
In other words, the spaces $\YY$ and $\YY/\mathrm{ker\ }T^\ast $
both specify the same constraint sets $\{\ell\in\LL; T\ell=x\}.$
\end{itemize}
The effective assumptions are the following ones.
\begin{itemize}
    \item[-] The specific form of the objective function $\Fs$ as a
convex conjugate makes it a convex $\sLU$-closed function.
    \item[-] \HFii\ and \HC\ are geometric restrictions.
    \item[-] \HTi\ is a
regularity assumption on $T.$
\end{itemize}

\subsection{Variants of $(P)$ and $(D)$}
\label{sec:problems} These variants are expressed below in terms
of new spaces and functions. Let us first introduce them.

\subsubsection*{The norms $|\cdot|_\F$ and $|\cdot|_\La$}
Let $\F_{\pm}(u)=\max(\F(u),\F(-u)).$ By \HFi\ and
\HFii, $\{u\in\UU; \F_{\pm}(u)\leq 1\}$  is a convex absorbing
balanced set. Hence its gauge functional which  is defined for all
$ u\in\UU$ by $|u|_\F \eqdef \inf\{\alpha>0;
\F_{\pm}(u/\alpha))\leq 1\}$ is a seminorm. Thanks to hypothesis
\HFiii, it is a norm.
\\
Taking \HTi\ into account, one can define
\begin{equation}
\label{La}
  \La(y)\eqdef \F(T^\ast y), y\in\YY.
\end{equation}
Let $\La_\pm(y)= \max(\La(y),\La(-y)).$ The gauge functional on
$\YY$ of the set $\{y\in\YY; \La_\pm(y)\leq 1\}$ is $|y|_\La
\eqdef \inf\{\alpha>0;\La_\pm(y/\alpha)\leq 1\}, y\in\YY.$ Thanks
to \HF\ and \HT, it is a norm and
\begin{equation}\label{eq-117}
    |y|_\La=|T^\ast y|_\F, \quad y\in\YY.
\end{equation}

\subsubsection*{The spaces}
Let
\begin{eqnarray*}
& &\Ui \mathrm{\ be\ the\ } \NF\textrm{-completion\ of\ } \UU\mathrm{\ and\ let }\\
  & &\Li\eqdef (\UU,\NF)'\textrm{\ be\ the\ topological\ dual\ space\ of\ }
  (\UU,\NF).
\end{eqnarray*}
Of course, we have $(\Ui,\NF)'\cong\Li\subset\LL$ where any $\ell$
 in $\Ui'$ is identified with its restriction to $\UU.$ Similarly, we introduce
\begin{eqnarray*}
& &\Yi \mathrm{\ the\ } \NL\textrm{-completion\ of\ } \YY
  \mathrm{\ and }\\
 & &\Xi\eqdef (\YY,\NL)'\mathrm{\ the\ topological\ dual\ space\ of\
 }(\YY,\NL).
\end{eqnarray*}
We have $(\Yi,\NL)'\cong\Xi\subset\XX$ where any $x$ in $\Yi'$ is
identified with its restriction to $\YY.$
\\
We also have to consider the algebraic dual space $\Uii$ and
$\Yii$ of $\Li$ and $\Xi.$

\subsubsection*{The adjoint operators of $T$} It will be proved at Lemma
\ref{L2} that
\begin{equation}\label{eq-101}
    T\Li\subset\Xi
\end{equation}
Let us denote $T_1$ the restriction of $T$ to $\Li\subset\LL.$ By
(\ref{eq-101}), we have $T_1: \Li\to\Xi.$ Let us define its
adjoint $T_2^\ast :\Yii\rightarrow\Uii$ for all $\omega\in\Yii$
by:
\begin{displaymath}
  \langle  \ell,T_2^\ast\omega\rangle_{\Li,\Uii}=\langle T_1\ell,\omega\rangle_{\Xi,\Yii},
  \forall \ell\in\Li.
\end{displaymath}
This definition is meaningful, thanks to (\ref{eq-101}). We denote
$T_1^*$ the restriction of $T_2^*$ to $\Yi\subset\Yii.$ Of course,
it is defined for any $y\in\Yi,$ by
$$
\langle \ell,T_1^\ast y\rangle_{\Li,\Li^\aldu}=\langle y,
T\ell\rangle_{\Yi,\Xi},\ \forall \ell\in\Li.
$$
It will proved at Lemma \ref{L2} that
\begin{equation}\label{eq-102}
    T^*_1\Yi\subset\Ui
\end{equation}
We have the inclusions $\YY\subset\Yi\subset\Yii.$ The adjoint
operators $T^*$ and $T^*_1$ are the restrictions of $T^*_2$ to
$\YY$ and $\Yi.$

\subsubsection*{Some modifications of $\F$ and $\La$}
The convex conjugate of $\F$  the dual pairing
$\langle\UU,\LL\rangle$ is
\begin{equation*}
    \Fs(\ell)\eqdef \sup_{u\in\UU}\{\ul-\F(u)\}, \ell\in\LL\\
\end{equation*}
We introduce the following modifications of $\F:$
\begin{eqnarray*}
    \F_0(u)&\eqdef&\sup_{\ell\in\LL}\{\ul -\Fs(\ell)\}, u\in\UU\\
    \F_1(u)&\eqdef&\sup_{\ell\in\Li}\{\ul-\Fs(\ell)\}, u\in\Ui\\
     \Fd(\zeta)&\eqdef&\sup_{\ell\in\Li}\{\lz-\Fs(\ell)\}, \zeta\in\Uii.
\end{eqnarray*}
They are respectively $\sUL,$ $\sULi$ and $\sULii$-closed convex
functions. It is immediate to see that the restriction of $\Fd$ to
$\Ui$ is $\F_1.$ As $\Li=\Ui',$ $\F_1$ is also the
$|\cdot|_{\F}$-closed convex regularization of $\F.$ The function
$\Fd$ is the extension $\bar\F$ which appears in the introductory
Section \ref{sec:introduction}.
\\
We also introduce
\begin{eqnarray*}
    \La_0(y)&\eqdef& \F_0(T^\ast y),\  y\in\YY\\
  \La_1(y)&\eqdef& \F_1(T_1^\ast y),\  y\in\Yi\\
   \Lb(\omega)&\eqdef&\Fd(T_2^\ast \omega),\ \omega\in\Yii
\end{eqnarray*}
which look like the definition (\ref{La}). Note that thanks to
\HTi\ and (\ref{eq-102}), the first equalities are meaningful.
Because of the previous remarks, the restriction of $\Lb$ to $\Yi$
is $\La_1.$

\subsubsection*{The optimization problems}
Let  $\Fs_0$ and $\F_1^*$ be the convex conjugates of $\F_0$ and
$\F_1$ with respect to the dual pairings $\langle\UU,\LL\rangle$
and $\langle\Ui,\Li\rangle:$
\begin{eqnarray*}
    \Fs_0(\ell)&\eqdef& \sup_{u\in\UU}\{\ul-\F_0(u)\}, \ell\in\LL\\
    \Fs_1(\ell)&\eqdef& \sup_{u\in\Ui}\{\ul-\F_1(u)\}, \ell\in\Li
\end{eqnarray*}
and $\Ls_0, \Ls_1$ be the convex conjugates of  $\La_0,\La_1$ with
respect to the dual pairings $\langle\YY,\XX\rangle$ and
$\langle\Yi,\Xi\rangle:$
\begin{eqnarray*}
     \Ls_0(x)&\eqdef& \sup_{y\in\YY}\{\yx-\La_0(y)\}, x\in\XX\\
    \Ls_1(x)&\eqdef& \sup_{y\in\Yi}\{\yx-\La_1(y)\}, x\in\Xi\\
\end{eqnarray*}
 Finally,  denote
 $$
  \Ci= C\cap\Xi.
$$
 The optimization problems to be considered are
\begin{align}
&\textsl{minimize } \Fs(\ell) & &\textsl{subject to } T\ell\in
C,\quad \ell\in\LL   \tag{$P$}\\
&\textsl{minimize } \Fs_1(\ell) & &\textsl{subject to }
T\ell\in\Ci,\quad \ell\in\Li   \tag{$P_1$}\\
&\textsl{minimize } \Ls_0(x) & &\textsl{subject to }
x\in\Ci,\quad x\in\Xi   \tag{$P_{1,\XX}$}\\
&\textsl{maximize } \inf_{x\in C}\yx - \La_0(y),  &
&y\in\YY  & &   \tag{$D_0$}\\
&\textsl{maximize } \inf_{x\in\Ci}\yx - \La_1(y),  &
&y\in\Yi  & &   \tag{$D_1$}\\
&\textsl{maximize } \inf_{x\in\Ci}\xo - \Lb(\omega),  &
&\omega\in\Yii & & \tag{$D_2$}
\end{align}

\subsection{Statement of the abstract results}

We are now ready to give answers to the questions related to $(P)$
and $(D)$ in an abstract setting.

\begin{theorem}[Primal attainment and dual equality]\label{xP3}
Assume that \HF\ and \HT\ hold.
\begin{enumerate}
    \item[(a)] For all $x$ in $\XX,$ we have the little dual equality
        \begin{equation}\label{xped}
            \inf\{\Fs(\ell); \ell\in\LL, T\ell=x\}=\Ls_0(x)\in
            [0,\infty].
        \end{equation}
        Moreover, in restriction to $\Xi,$ $\Ls_0=\Ls_1$ and
        $\Ls_1$ is $\sigma(\Xi,\Yi)$-inf-compact.
    \item[(b)] The problems \Po\ and \Pi\ are equivalent: they have
        the same solutions and $\inf(P)=\inf(P_1)\in[0,\infty].$
    \item[(c)] If $ C$ is convex and $\sXY$-closed, we have the dual equality
\begin{displaymath}
  \inf(P)=\sup(D_0)\in[0,\infty].
\end{displaymath}
\end{enumerate}
Assume that \HF, \HT\ and \HC\ hold.
\begin{enumerate}
    \item[(d)]  We have the dual equalities
\begin{equation}
   \inf(P)=\inf(P_1)=\sup(D_1)=\sup(D_2)=\inf_{x\in C}\Ls_0(x)\in [0,\infty]\label{xed1}
\end{equation}
    \item[(e)] If in addition $\inf(P)<\infty,$ then \Po\ is attained in
    $\Li.$ Moreover, any minimizing sequence for \Po\ has $\sLUi$-cluster
    points and every such cluster point solves \Po.
    \item[(f)] Let $\lb\in\Li$ be a solution to $(P),$ then
$\xb\eqdef T\lb$ is a solution to \PX\ and
$\inf(P)=\Fs(\lb)=\Ls_0(\xb).$
\end{enumerate}
\end{theorem}

\begin{theorem}[Dual attainment and representation. Interior convex constraint]\label{T3a}
Assume that \HF, \HT\ and \HC\ hold and also suppose that the
interior constraint qualification
 \begin{equation}\label{xeq-106}
    C\cap\icordom\Ls_0\not=\emptyset
\end{equation}
is satisfied. Then, the following statements hold true.
\begin{enumerate}[(a)]
    \item  The primal problem \Po\ is attained in $\Li$ and the dual problem \Dii\ is attained in $\Yii$
    \item Any  $\lb\in\Li$ is a solution to \Po\
 if and only if there exists $\ob\in\Yii$ such that the following three statements hold
 \begin{equation*}
    \left\{\begin{array}{cl}
      (1) & T\lb\in C \\
      (2) & \langle T\lb,\ob\rangle\leq \langle x,\ob\rangle
      \textrm{ for all }x\in C_1 \\
      (3) & \lb\in\partial_{\Li}\Fd(T^*_2\ob) \\
    \end{array}\right.
\end{equation*}
More, these three statements hold if and only if: $\lb$ is a
solution to \Po, $\ob$ is a solution to \Dii\ and
$\inf(P)$=$\sup(D_2).$
\\ It is
well-known that the representation formula
        \begin{equation}\label{xeq-110}
    \lb\in\partial_{\Li}\Fd(T^*_2\ob)
\end{equation}
 is equivalent to Young's identity
\begin{equation}\label{xeq-107}
    \Fs(\lb)+\Fd(T^*_2\ob)=\langle T\lb,\ob\rangle.
\end{equation}
    \item Any solution $\ob$ of \Dii\ shares the following properties
\begin{enumerate}[(1)]
 \item $\ob$ stands in the $\sYXx$-closure of $\dom\La_1.$
 \item $T_2^\ast \ob$ stands in the $\sULii$-closures of
 $T_1^\ast(\dom\La_1)$ and $\dom\F.$
 \item For any $x_o$ in $C\cap\icordom\Ls_1,$ $\ob$ is $j_{D_{x_o}}$-\usc\     and $j_{-D_{x_o}}$-lower semicontinuous at $0,$ where   $j_{D_{x_o}}$ and
    $j_{-D_{x_o}}$ are the gauge functionals on $\Xi$ of the convex sets $D_{x_o}=\{x\in\Xi;\Ls_0(x_o+x)\leq\Ls_0(x_o)+1\}$ and
    $-D_{x_o}$.
        \end{enumerate}
\end{enumerate}
\end{theorem}

As will be seen at Section \ref{sec:transport}, the
Monge-Kantorovich problem provides an important example where no
constraint is interior (see Remark \ref{rem-04}). In order to
solve it without imposing constraint qualification, we are going
to consider the more general situation (\ref{eq-14b}) where the
constraint is said to be a \emph{subgradient constraint}. This
means that $ \xb\in\diffdom\Ls_0$ with
\begin{eqnarray*}
    \diffdom\Ls_0&=&\{x\in\Xi; \partial_{\Yii}\Ls_0(x)\not
    =\emptyset\}\quad \textrm{ where}\\
    \partial_{\Yii}\Ls_0(x)&=&\{\omega\in\Yii; \Ls_0(x')\geq\Ls_0(x)+\langle x'-x,\omega\rangle,\forall
    x'\in\Xi\}.
\end{eqnarray*}

 Two new optimization problems to be considered are
\begin{align}
&\textsl{minimize } \Fs(\ell) & &\textsl{subject to } T\ell=\xb,\quad \ell\in\LL   \tag{$P^{\xb}$}\\
&\textsl{maximize } \langle\xb,\omega\rangle-\La_2(\omega), &
&\omega\in\Yii & & \tag{$D_2^{\xb}$}
\end{align}
where $\xb\in\XX.$ This corresponds to the simplified case where
$C$ is reduced to the single point $\xb.$

\begin{theorem}[Dual attainment and representation. Subgradient affine constraint]\label{T3b}
Let us assume that \HF\ and \HT\ hold and suppose that
$\xb\in\dom\Ls_0.$ Then, $\inf(P^{\xb})<\infty.$  If in addition,
\begin{equation}\label{xCQbis}
    \xb\in\diffdom\Ls_0,
\end{equation}
then the following statements hold true.
\begin{enumerate}[(a)]
    \item The primal problem $(P^{\xb})$ is attained in $\Li$ and the dual problem \Dxb\ is attained in $\Yii.$
     \item Any  $\lb\in\Li$ is a solution to $(P^{\xb})$ if and only
     if $T\lb=\xb$ and
     there exists $\ob\in\Yii$ such that (\ref{xeq-110})
     or equivalently (\ref{xeq-107}) holds.
     \\
     More, this occurs if and only if: $\lb$ is a
     solution to \Po, $\ob$ is a solution to \Dxb\ with $\xb:=T\lb$ and
    $\inf(P^{\xb})=\sup(D_2^{\xb}).$
     \item Any solution $\ob$ of \Dxb,   shares the following properties
\begin{enumerate}[(1)]
 \item $\ob$ stands in the $\sYXx$-closure of $\dom\La_1.$
 \item $T_2^\ast \ob$ stands in the $\sULii$-closures of
 $T_1^\ast(\dom\La_1)$ and $\dom\F.$
 \item Let $\ob$ be any solution of \Dxb\ with $\xb\in\icordom\Ls_0$. Then, $\ob$ is $j_{D_{\xb}}$-\usc\             and $j_{-D_{\xb}}$-lower semicontinuous at $0$ where $D_{\xb}=\{x\in\Xi;\Ls_0(\xb+x)\leq\Ls_0(\xb)+1\}.$
        \end{enumerate}
 \end{enumerate}
\end{theorem}

\section{Application to the Monge-Kantorovich optimal transport problem}\label{sec:transport}

We apply the results of Section \ref{sec:abstractpb} to the
Monge-Kantorovich problem. Recall that we take $A$ and $B$ two
Polish spaces furnished with their  Borel $\sigma$-fields. Their
product space $\AB$ is endowed with the product topology and the
corresponding Borel $\sigma$-field.  The \lsc\ cost function
$c:\AB \to [0,\infty]$ may take infinite values. Let us also take
two probability measures $\mu\in\PA$ and $\nu\in\PB$ on $A$ and
$B.$ The Monge-Kantorovich problem is
\begin{equation*}
    \textsl{minimize } \pi\in\PAB\mapsto \IAB c(a,b)\,\pi(dadb)
    \textsl{ subject to }\pi\in P(\mu,\nu) \tag{$\MK$}
\end{equation*}
where $P(\mu,\nu)$ is the set of all $\pi\in\PAB$ with prescribed
marginals $\pi_A=\mu$ on $A$ and $\pi_B=\nu$ on $B.$

\subsection{Statement of the results}

Let us fix some notations.
 We denote $C_A,$ $C_B$ and $\CAB$ the spaces of all continuous
bounded functions on $A,$ $B$ and $\AB.$ The  Kantorovich
maximization problem:
\begin{equation}\label{eq-28}
\begin{split}
   &\textsl{maximize }\int_A\varphi\,d\mu+\int_B\psi\,d\nu \textsl{ for all
    }\varphi,\psi \textsl{ such that }\\
    &\varphi\in C_A,\psi\in C_B\textsl{ and }\varphi\oplus\psi\leq
    c\\
\end{split}
\tag{$K$}
\end{equation}
is the basic dual problem of $(\MK).$  We also consider the
following extended version of $(K):$
\begin{equation}\label{eq-25}
   \begin{split}
     &\textsl{maximize }\int_{A}\varphi\,d\mu+\int_{B}\psi\,d\nu \textsl{ for all
    }\varphi\in \mathbb{R}^{A},\psi\in\mathbb{R}^{B} \textsl{ such that }\\
     &\varphi\in L_1(A,\mu), \psi\in L_1(B,\nu)
     \textrm{ and }\varphi\oplus\psi\leq c\textrm{ everywhere on
     }\AB.
    \end{split}\tag{$\overline{K}$}
    \end{equation}

\begin{remark} The real-valued function $\varphi\in \mathbb{R}^{A}$ is defined
    \emph{everywhere}, rather than  $\mu$-almost everywhere, and $\varphi\in L_1(A,\mu)$ implies that it is \emph{$\mu$-measurable}.
    This means that there exists
    some measurable set $N_A\subset A$ such that $\mu(N_A)=0$ and
    $\1_{N_A}\varphi$ is measurable. A similar remark holds for $\psi.$
\end{remark}

The set of all probability measures $\pi$ on $\AB$ such that $\IAB
c\,d\pi<\infty$ is denoted $\Pc.$ By Definition \ref{def-02}, an
optimal plan stands in $\Pc.$ In the next theorem, $\Pc$ will be
endowed with the weak topology $\sigma(\Pc,\mathcal{C}_c)$ where
$\mathcal{C}_c$ is the space of all continuous functions $u$ on
$\AB$ such that $|u|\leq k(1+c)$ for some $k\geq 0.$

\begin{theorem}[Dual equality and
primal attainment]\label{res-MK1}\
\begin{enumerate}
 \item The dual equality for $(\MK)$ is
$$
    \inf(\MK)=\sup(K)=\sup(\overline{K})\in [0,\infty].
$$
 \item Assume that there exists
some $\pi^o$ in $P(\mu,\nu)$ such that $\IAB c\,d\pi^o<\infty.$
Then:
    \begin{enumerate}
            \item There is at least an optimal plan and all the optimal plans are in $\Pc;$
            \item Any minimizing sequence is relatively compact for the topology
            $\sigma(\Pc, \mathcal{C}_c)$ and all its cluster points are optimal plans.
   \end{enumerate}
\end{enumerate}
\end{theorem}

This result is well-known. The dual equality
$\inf(\MK)=\sup(K)=\sup(\overline{K})$ is the Kantorovich dual
equality. The proof of Theorem \ref{res-MK1} will be an
opportunity to make precise the abstract material $\F,$ $\UU,$
$T\dots$  in terms of the Monge-Kantorovich problem.

Next, we state the characterization of the optimal plans without
restriction.

\begin{theorem}[Characterization of the optimal
plans]\label{res-MK2}\
\begin{enumerate}
    \item A probability measure $\pi\in\PAB$ is an optimal plan if and only if there exist
     two finitely-valued functions $\varphi\in\mathbb{R}^A$ and $\psi\in\mathbb{R}^B$  such that
    \begin{equation}\label{eq-26}
    \left\{\begin{array}{ll}
   (a) &\pi_A=\mu, \pi_B=\nu, \IAB c\,d\pi<\infty,\\
   (b)   &\varphi\oplus\psi\leq c \quad \textrm{everywhere and}\\
    (c)   &\varphi\oplus\psi= c \quad \pi\textrm{-almost everywhere.} \\
    \end{array}\right.
    \end{equation}
\item Let $\varphi$ and $\psi$ be finitely-valued functions on $A$
and $B$ and let be $\pi\in\PAB.$
\begin{enumerate}[(a)]
    \item  If $\varphi$ is $\mu$-measurable and $\psi$ is
    $\nu$-measurable, the following statements are equivalent:
        \begin{itemize}
            \item[-] $\varphi,$ $\psi$ and $\pi$ satisfy
            (\ref{eq-26});
             \item[-] $\pi$ is an optimal plan and $(\varphi,\psi)$ is a solution of $(\overline{K}).$
        \end{itemize}
    \item In the general case where $\varphi$ and $\psi$ are not assumed
    to be
measurable, consider the following statements:
\begin{enumerate}[(i)]
            \item $\varphi,$ $\psi$ and $\pi$ satisfy
            (\ref{eq-26});
             \item $\ls\varphi,$ $\ls\psi$ and $\pi$ satisfy
            (\ref{eq-26});
             \item $\pi$ is an optimal plan and $(\ls\varphi,\ls\psi)$ is a solution of $(\overline{K}).$
        \end{enumerate}
        Then: $(i)\Rightarrow(ii)\Leftrightarrow(iii).$
\end{enumerate}
\end{enumerate}
\end{theorem}
This new result improves the already existing literature on the
subject. It is important to note that \emph{the functions
$\varphi$ and $\psi$ satisfying (\ref{eq-26}) are neither assumed
to be integrable nor to be measurable}. Next theorem shows that
they can be further specified.
\begin{theorem}[More about necessary conditions]\label{res-MK3}
Assume that  $\pi$ is an optimal plan. Then, there exist two
finitely-valued functions $\varphi\in\mathbb{R}^A$ and
$\psi\in\mathbb{R}^B$  such that $\varphi\in L_1(A,\mu),$ $\psi\in
L_1(B,\nu)$ and
\begin{equation*}
    \left\{\begin{array}{clll}
      |\varphi\oplus\psi|&\leq &c & \textrm{everywhere and} \\
       \varphi\oplus\psi&= &c &  \textrm{on }\suppp\bigcap\{c<\infty\}.\\
    \end{array}\right.
\end{equation*}
Clearly, $(\varphi,\psi)$ is a maximizer of $(\overline{K}).$
\end{theorem}

\begin{remarks}\
\begin{enumerate}[a.]
    \item Note that  any optimal plan $\pi$ is satisfies $\suppp\subset\cl\{c<\infty\}.$
    \item Recall that $\pi$ is said be to \emph{concentrated} on the
measurable set $\Gamma$ if $\pi(\Gamma)=1.$ For instance,
(\ref{eq-26}-c) is equivalent to the existence of some set
$\Gamma$ on which $\pi$ is concentrated and $\varphi\oplus\psi=c$
on $\Gamma.$ The \emph{support} of $\pi,$ denoted $\suppp,$ is the
closure of the union of all the sets $\Gamma$ on which $\pi$ is
concentrated.
\end{enumerate}
\end{remarks}

\subsection{Proof of Theorem \ref{res-MK1}}
We apply the general results of Section \ref{sec:abstractpb}.

\subsubsection*{The operators $T$ and $T^*$}

The algebraic dual spaces of $C_A,$ $C_B$ and $\CAB$ are $\CAs,$
$\CBs$ and $\CABs.$ We define the marginal operator
\begin{equation*}
    T\ell=(\ell_A,\ell_B)\in\CAs\times\CBs,\quad \ell\in\CABs
\end{equation*}
where $\langle\varphi,\ell_A\rangle=\langle\varphi\otimes
1,\ell\rangle$ and $\langle\psi,\ell_B\rangle=\langle
1\otimes\psi,\ell\rangle$ for all $\varphi\in C_A$ and all
$\psi\in C_B.$
\\
Let us identify the operator $T^\ast.$ For all $(\varphi,\psi)\in
C_A\times C_B$ and all $\ell\in \LL,$ we have $\langle
T^\ast(\varphi,\psi),\ell\rangle_{\LL^*,\LL}=\langle
(\varphi,\psi),(\ell_A,\ell_B)\rangle
 =\langle\varphi,\ell_A\rangle+\langle\psi,\ell_B\rangle
 =\langle\varphi\oplus\psi,\ell\rangle_{\UU,\LL}$ where $\varphi\oplus\psi(a,b)=\varphi(a)+\psi(b).$
Hence, for each $\varphi\in C_A$ and $\psi\in C_B,$
\begin{equation}\label{eq-16}
    T^\ast(\varphi,\psi)=\varphi\oplus\psi\in C_{A\times B}.
\end{equation}

\subsubsection*{The problem $(P)$}

 Then, the Diagram 0 is built with $\UU=\CAB,$
$\LL=\CABs,$ $\XX=\CAs\times\CBs$ and $\YY=C_A\times C_B.$ Here
and below, we denote the convex indicator function of the set $X,$
$$\delta_X(x)=\left\{\begin{array}{ll}
  0 & \textrm{if } x \in X \\
  +\infty & \textrm{otherwise.} \\
\end{array}\right.$$
 Choosing $C=\{(\mu,\nu)\}$ and
 $ 
    \Phi(u)=\delta_{\{u\leq c\}}
    , u\in\CAB
 $ 
 we get
 $ 
    \Fs(\ell)=\sup\{\langle u,\ell\rangle; u\in\CAB,u\leq
    c\}, \ell\in\CABs
 $ 
 and we obtain the primal problem
\begin{equation}
    \textsl{minimize } \Fs(\ell) \textsl{ subject to }
    \ell_A=\mu \textsl{ and } \ell_B=\nu
    ,\quad \ell\in\CABs.   \tag{$P$}
\end{equation}
It will be shown at Proposition \ref{res-02} that the
corresponding problem \Pi\ is $(\MK).$

\subsubsection*{The problem $(D_0)$}

Now, let's have a look at $\Phi_0.$ As $\{u\in\CAB; u\leq c\}$ is
convex and $\sigma(\CAB,\CABs)$-closed,
 we have $\Phi_0=\Phi.$ Therefore,  for each $\varphi\in C_A$ and $\psi\in C_B,$
\begin{equation*}
    \La_0(\varphi,\psi)=\La(\varphi,\psi)=\Phi(T^\ast(\varphi,\psi))=\delta_{\{\varphi\oplus\psi\leq c\}}
\end{equation*}
and the dual problem is
\begin{equation}
    \textsl{maximize } \int_A\varphi\,d\mu + \int_B\psi\,d\nu \textsl{ subject to }
    \varphi\oplus\psi\leq c
    ,\quad \varphi\in C_A,\psi\in C_B  \tag{$D_0$}
\end{equation}
whose value is
\begin{equation*}
    \Ls(\mu,\nu)=\sup\left\{\int_A\varphi\,d\mu + \int_B\psi\,d\nu; \varphi\in C_A,\psi\in
    C_B :
    \varphi\oplus\psi\leq c\right\}
\end{equation*}
As $\La=\La_0$ and $\Ls_0=\Ls_1$ (Theorem \ref{xP3}-a), we have:
$\Ls_0=\Ls_1=\Ls.$

\subsubsection*{The hypotheses $(H)$}

We begin with a simple remark.
\begin{remark}\label{rem-03}
One can choose $c\geq 1$ without loss of generality. Indeed, with
$c\geq 0$ taking $\tilde{c}=c+1$ one obtains $\IAB c\,d\pi=\IAB
\tilde{c}\,d\pi-1$ for all $\pi\in\PAB.$ Consequently, the
minimization problems $(\MK)$ and $(\widetilde{\MK})$ associated
with $c$ and $\tilde{c}$ share the same minimizers and their
values are related by $\inf(\MK)=\inf(\widetilde{\MK})-1.$ \\ It
also follows from these considerations that our results still hold
under the assumption that $c$ is bounded below rather than $c$ is
\positif.
\end{remark}
 We assume
from now on that $c\geq 1.$ This guarantees \HFii. In the case
where $c$ is finitely valued, the remaining hypotheses $(H)$
follow by (\ref{eq-16}) and direct inspection.

If $c$ is infinite somewhere, then \HFiii\ fails. Indeed, for any
function $u\in\CAB,$ we have $\F(tu)=0$ for all real $t$ if and
only if $\{u\not=0\}\subset\AB\setminus\s$ where
\begin{equation*}
    \s=\cl \{c<\infty\}
\end{equation*}
is the closure of $\{(a,b)\in\AB; c(a,b)<\infty\}.$ The way to get
rid of this problem is standard. Let $u\sim v$ be the equivalence
relation on $\mathbb{R}^{\AB}$ defined by $\uS=v_{|\s},$ i.e. $u$
and $v$ match on $\s.$ The space $\UU$ to be considered is the
factor space
$$\UU:=\CAB/\sim$$ Clearly, if $u\sim v$ then $\F(u)=\F(v).$
Hence, it is possible to identify without loss of generality any
$u\in\CAB$ with its equivalence class which in turn is identified
with the restriction $\uS$ of $u$ to $\s.$

\subsubsection*{The problem $(P_1)$}
Recall that $c\geq 1$ without loss of generality. Let us first
identify the space $\Li.$ As $\Phi_\pm(u)=\delta_{\{|u|\leq c\}},$
we obtain the seminorm $|u|_\Phi=\sup|u/c|:=\|u\|_c$ on $\CAB$
which becomes a norm on $\UU,$
\begin{eqnarray*}
 \Ui&=&\Cc:=\{\uS; u:\AB\to\mathbb{R}, u \textrm{ continuous and }|u|\leq
kc \textrm{ for some real }k\} \\
  \Li&=&(\Cc,\|\cdot\|_c)'=\Cc'.
\end{eqnarray*}
Obviously, any $\pi$ in $\Pc$ has its support included in $\s$ and
belongs to $\Cc'$ with the dual bracket $\langle
\uS,\pi\rangle=\IS u\,d\pi,$ $\uS\in\Cc.$ In what follows, it will
be written equivalently
\begin{itemize}
    \item $\uS\in\Cc$ to
specify that the equivalence class of $u$ stands in $\Cc$ and
    \item $u\in\Cc$ to specify that the restriction $\uS$ of the
    continuous function $u$ on $\AB$ stands in $\Cc.$
\end{itemize}
Clearly, the function $\F_1$ is
\begin{equation*}
    \F_1(u)=\delta_{\{u\leq c\}},\quad u\in \Cc
\end{equation*}
and the modified primal problem is
\begin{equation}
    \textsl{minimize } \Fs_1(\ell) \textsl{ subject to
    }\ell_A=\mu \textsl{ and } \ell_B=\nu,\quad \ell\in \Cc'
    \tag{$P_1$}
\end{equation}
where  for each $\ell\in \Cc',$
\begin{equation*}
  \Fs_1(\ell)=\sup\{\langle u,\ell\rangle; u\in \Cc, u\leq
    c\}.
\end{equation*}

\begin{remark}\label{rem-02}
Two representations of $\Cc'.$
\begin{enumerate}[a.]
    \item Let  $\Cc(\s)$ be the space of all continuous functions $w$ on
$\s$ (\hbox{w.r.t.} the relative topology) such that
$\|w\|_c=\sup_{\s}|w/c|<\infty$ and $\Cc(\s)'$ be the topological
dual space of the normed space $(\Cc(\s),\|\cdot\|_c).$ Let
$\mathcal{E}$ be the subspace of all functions in $\Cc(\s)$ which
can be continuously extended to the whole space $\AB.$ There is a
one-one correspondence between $\Cc'$ and the dual space
$(\mathcal{E},\|\cdot\|_c)'.$
    \item There is also a one-one correspondence between $\Cc'$ and the space
    of all linear forms $\ell$ on the space of all continuous functions on $\AB$
    such that $\supp\ell\subset\s$ (see Definition \ref{def-01} below) and
    $\sup\{\langle u,\ell\rangle=\langle \uS,\ell\rangle; u:\|u\|_c\leq
    1\}<\infty.$
\end{enumerate}
\end{remark}

\begin{definition}\label{def-01}
For any linear form $\ell$ on the space of all continuous
functions on $\AB,$ we define the \emph{support} of $\ell$ as the
subset of all $(a,b)\in\AB$ such that for any neighborhood $G$ of
$(a,b),$ there exists some function $u$ in $\CAB$ satisfying
$\{u\not=0\}\subset G$ and $\langle u,\ell\rangle \not=0.$ It is
denoted $\mathrm{supp\,}\ell.$
\end{definition}

\begin{definition}\
\begin{enumerate}[(a)]
    \item One says that $\ell\in \Cc'$ \emph{acts as a probability measure}
if there exists $\tilde{\ell}\in\PAB$ such that
$\supp\tilde{\ell}\subset\s$ and for all $u\in \CAB,$ $\langle
\uS,\ell\rangle=\IS u\,d\tilde{\ell}.$  In this case, we write:
$\ell\in\PS.$
    \item One says that $\ell\in \Cc'$ stands in $\Pc$ if there exists
$\tilde{\ell}\in\Pc$ such that  for all $u\in \Cc,$ $\langle
\uS,\ell\rangle=\IS u\,d\tilde{\ell}.$  In this case, we write:
$\ell\in\Pc.$
\end{enumerate}
\end{definition}
Of course, if there exists $\tilde{\ell}$ satisfying (a), it
belongs to $\Pc$ and is unique since any probability measure on a
metric space is determined by its values on the continuous bounded
functions. This explains why the notation $\ell\in\Pc$ in (b)
isn't misleading.
\\
Note also that any probability measure $\tilde{\ell}\in\Pc$ has a
support included in $\s.$ Since $\AB$ is a metric space, for any
$\ell\in\Pc$ acting as a measure, $\supp\ell$ in the sense of
Definition \ref{def-01} matches with the usual support of the
measure $\tilde{\ell}.$

\subsubsection*{Completing the proof of Theorem \ref{res-MK1}}
The full connection with the Monge-Kantorovich problem is given by
the following Proposition \ref{res-02}. Clearly, with this
proposition in hand, Theorem \ref{res-MK1} directly follows from
Theorem \ref{xP3} and the obvious inequalities
$\sup(K)\leq\sup(\overline{K})\leq\inf(\MK).$
\begin{proposition}\label{res-02}
For all $\ell\in \Cc',$
\begin{listalpha}
    \item $\Fs_1(\ell)<\infty\Rightarrow \ell\geq 0,$
    \item $\Fs_1(\ell)<\infty\Rightarrow \supp\ell\subset\s,$
    \item $[\ell\geq 0,\supp\ell\subset\s, \ell_A=\mu \textsl{ and } \ell_B=\nu]\Rightarrow
    \ell\in\PS$ and
    \item for all $\ell\in\PS,$ $\Fs_1(\ell)=\IS c\,d\ell.$
\end{listalpha}
It follows that \begin{itemize}
    \item[-] $\dom\Fs_1\subset\Pc$ and
    \item[-] the problems $(\MK)$ and $(P_1)$
share the same values and the same minimizers.
\end{itemize}
\end{proposition}

 \proof
 Clearly, the last statement follows from the first part of the
 proposition. The proof is divided into four parts.

\Boulette{a} Suppose that $\ell\in \Cc'$ isn't in the \positif\
cone. This means that there exists $u_o\in \Cc$ such that $u_o\geq
0$ and $\langle u_o,\ell\rangle<0.$ Since  $u_o$ satisfies
$\lambda u_o\leq 0\leq c$ for all $\lambda<0,$ we have
$\Fs_1(\ell)\geq \sup_{\lambda<0}\{\langle \lambda
u_o,\ell\rangle\}=+\infty.$ Hence, $\Fs_1(\ell)<\infty$ implies
that $\ell\geq 0$ and one can restrict our attention to the
\positif\  $\ell$'s.

\Boulette{b} Suppose ad absurdum that $\supp\ell\varsubsetneq\s.$
Then, there exists a \positif\ function $u_o\in\CAB$ such that
$\{u_o>0\}\cap\s = \emptyset$ and $\langle u_o,\ell\rangle>0.$ As
$tu_o\leq c_{|\AB\setminus\s}\equiv\infty$ for all $t>0,$
$\Fs_1(\ell)\geq \sup_{t>0}\{\langle tu_o,\ell\rangle\}=+\infty.$

\Boulette{c} Let us take $\ell\geq 0$ such that
$\supp\ell\subset\s,$ $\ell_A=\mu$ and $\ell_B=\nu.$ It is clear
that $\langle 1,\ell\rangle=1.$ It remains to check that for any
$\ell\in \Cc'$
\begin{equation}\label{eq-18}
    [\ell\geq 0, \supp\ell\subset\s, \ell_A=\mu \textsl{ and } \ell_B=\nu]\Rightarrow
    \ell \textrm{ is }\sigma\textrm{-additive,}
\end{equation}
rather than only additive. Since $\AB$ is a metric space, one can
apply an extension of the construction of Daniell's integrals
(\cite{Neveu}, Proposition II.7.2) to see that $\ell$ acts as a
measure if and only if for any \decroissant\  sequence $(u_n)$ of
continuous functions such that $0\leq u_n\leq 1$ for all $n$ and
$\lim_{n\rightarrow\infty} u_n=0$ pointwise, we have
$\lim_{n\rightarrow\infty}\langle {u_n},\ell\rangle=0.$ This
insures the $\sigma$-additivity of $\ell.$ Note that as
$\supp\ell\subset\s,$ for all $u\in\Cc$ one can shortly write
$\langle u,\ell\rangle$ instead of the meaningful bracket $\langle
\uS,\ell\rangle.$
\\
Unfortunately, this pointwise convergence of $(u_n)$ is weaker
than the uniform convergence with respect to which any $\ell\in
\Cc'$ is continuous. Except if $\AB$ is compact, since in this
special case, any \decroissant\  sequence of continuous functions
which converges pointwise to zero also converges uniformly on the
compact space $\s.$
\\
So far, we have only used the fact that $\AB$ is a metric space.
We now rely on the Polishness of $A$ and $B$ to get rid of this
compactness restriction. It is known that any probability measure
$P$ on a Polish space $X$ is tight (i.e. a Radon measure): for all
$\epsilon>0,$ there exists a compact set $K_\epsilon\subset X$
such that $P(X\setminus K_\epsilon)\leq \epsilon$ (\cite{Neveu},
Proposition II.7.3). As in addition a Polish space is completely
regular, there exists a continuous function $f_\epsilon$ with a
compact support such  that $0\leq f_\epsilon\leq 1$ and $\int_X
(1-f_\epsilon)\,dP\leq \epsilon.$ This is true in particular for
the probability measures $\mu\in\PA$ and $\nu\in\PB$ which specify
the constraint in $(\MK).$ Hence, there exist $\varphi_\epsilon\in
C_A$ and $\psi_\epsilon\in C_B$ with compact supports such that
$0\leq \varphi_\epsilon,\psi_\epsilon\leq 1$ and $0\leq\int_A
(1-\varphi_\epsilon)\,d\mu,\int_B (1-\psi_\epsilon)\,d\nu \leq
\epsilon.$ It follows that any $\ell\in \Cc'$ with $\ell_A=\mu$
and $\ell_B=\nu$ satisfies $0\leq \langle
(1-\varphi_\epsilon\otimes\psi_\epsilon),\ell\rangle\leq
2\epsilon.$ With the following easy estimate $0\leq \langle
u_n,\ell\rangle\leq 2\epsilon +\langle
u_n(\varphi_\epsilon\otimes\psi_\epsilon),\ell\rangle $ and the
compactness of the support of
$\varphi_\epsilon\otimes\psi_\epsilon,$ one concludes that
$\lim_{n\rightarrow\infty}\langle u_n,\ell\rangle=0$ which proves
(\ref{eq-18}).

\Boulette{d} As $c$ is bounded below and \lsc\ on a metric space,
it is the pointwise limit of a \croissant\  sequence $(c_n)$ of
continuous bounded functions. It follows from the monotone
convergence theorem that for any $\ell\in\PS,$ $\Fs_1(\ell)=\IS
c\,d\ell.$ This completes the proof of the proposition.
\endproof

\subsection*{Optimal plan: an overview of the proofs of Theorems \ref{res-MK2} and \ref{res-MK3}}

The proofs of these theorems are postponed to Section
\ref{sec-04}. We first derive preliminary results at Sections
\ref{sec-01}, \ref{sec-02} and \ref{sec-03}.

At Section \ref{sec-01}, the abstract results of Section
\ref{sec:abstractpb} are translated in terms of the
Monge-Kantorovich problem. This is summarized at Theorem
\ref{res-20} which states an abstract characterization of the
optimal plans. This theorem directly results from the extended
saddle-point method. In particular, the optimal plan $\pi$ is
related to some linear form $\omega\in\Yii.$ It remains to show
that $\omega$ is the extension of some couple of functions
$(\varphi,\psi).$ This is done at Section \ref{sec-02} for the
sufficient condition and at Section \ref{sec-03} for the necessary
condition. The main results of Sections \ref{sec-02} and
\ref{sec-03} are respectively Lemma \ref{res-13} and Lemma
\ref{res-11c}.

\subsection{Optimal plan: applying the extended saddle-point
method}\label{sec-01}

The main result of this section is Theorem \ref{res-20} which
gives an abstract characterization of an optimal plan.

\subsubsection*{The space $\Xi$}

By (\ref{eq-117}), we see that
$|(\varphi,\psi)|_\Lambda=\|\varphi\oplus\psi\|_c.$ This leads to
\begin{equation*}
    \Xi=\{(\kappa_1,\kappa_2)\in C_A^*\times C_B^*; |(\kappa_1,\kappa_2)|_{\Lambda}^*<\infty\}
\end{equation*}
where $
  |(\kappa_1,\kappa_2)|_{\Lambda}^*
=
\sup\{\langle\varphi,\kappa_1\rangle+\langle\psi,\kappa_2\rangle;
(\varphi,\psi)\in\Yi, \|\varphi\oplus\psi\|_c\leq 1\}.$ The dual
equality (\ref{xed1}) gives
\begin{equation*}
    |(\kappa_1,\kappa_2)|_{\Lambda}^*
  =  \inf\left\{\left\|\ell\right\|_{c}^*; \ell\in \Cc': \ell_A=\kappa_1,\ell_B=\kappa_2\right\}.
\end{equation*}
Note that $\Xi$ is the space of all $(\kappa_1,\kappa_2)\in
C_A^*\times C_B^*$ such that $\kappa_1=\ell_A$ and
$\kappa_2=\ell_B$ for some $\ell$ in $\Cc'.$ Recall that the
elements of $\Xi=\Yi'$ are identified with their restriction to
$\YY$ which is dense in $\Yi.$

\begin{remark}[The space $\Yi$ and the problem \Di] The exact
description  of $\Yi$ and \Di\  will not be used later.
Nevertheless, as an illustration of our general results, we
describe them assuming that $c$ is finitely valued.
 As $\Lambda(\varphi,\psi)=\delta_{\{\varphi\oplus\psi\leq c\}}$,
one sees that
\begin{equation*}
    \Yi=\{(\varphi,\psi); \varphi:A\to\mathbb{R} \textrm{ continuous},
    \psi:B\to\mathbb{R} \textrm{ continuous}: \varphi\oplus\psi\in \Cc\}
\end{equation*}
This result is not as obvious as it seems to be. It follows from
an interesting paper \cite{BL92} of J.M.~Borwein and A.S.~Lewis
which studies the convergence of sequences of the form
$(\varphi_n\oplus\psi_n)_{n\geq 1}.$ The additive form
$\varphi\oplus\psi$ in the expression of $\Yi$ is proved at
(\cite{BL92}, Corollary 3.5) and the continuity of $\varphi$ and
$\psi$ is a consequence of (\cite{BL92}, Proposition 5.1).

The corresponding problem \Di\ is
\begin{equation*}
   \begin{split}
    &\textsl{maximize }\int_A\varphi\,d\mu+\int_B\psi\,d\nu \textsl{ for all
    }\varphi,\psi \textsl{ such that }\\
    &\varphi, \psi \textsl{ continuous, }\varphi\oplus\psi\in \Cc\textsl{ and }\varphi\oplus\psi\leq
    c.
    \end{split}
\end{equation*}
Anyway, we won't use this dual problem since it is sandwiched
between \Do\ and \Dii.
\end{remark}

\subsubsection*{The extension $\Fd$}
To proceed, one has to compute the extension $\Fd.$ As it is the
greatest convex $\sigma(\Cc^{\prime *},\Cc')$-\lsc\ extension of
$\Phi,$ we have
\begin{equation}\label{eq-72}
    \Fd(\xi)=\delta_{\overline{\Gamma}}(\xi),\quad \xi\in \Cc^{\prime *}
\end{equation}
where $\overline{\Gamma}$ is the $\sigma(\Cc^{\prime
*},\Cc')$-closure of
\begin{equation*}
    \Gamma=\{u\leq c\}\subset \Cc.
\end{equation*}
Any $\omega\in\Yii$ is decomposed as $\omega=(\omega_A,\omega_B)$
where for all $(\kappa_1,\kappa_2)\in\Xi,$
 $\langle\omega,(\kappa_1,\kappa_2)\rangle
 =\langle\bar{\omega},(\kappa_1,0)\rangle+\langle\bar{\omega},(0,\kappa_2)\rangle
 =\langle\omega_A,\kappa_1\rangle
+\langle\omega_B,\kappa_2\rangle$ where $\omega\in\Yii$ is seen as
the restriction to $\Xi$ of some linear form $\bar{\omega}$ on
$\XX=C_A^*\times C_B^*.$ The adjoint operator $T_2^*$
 is defined for all $\omega\in\Yii$ and $\ell\in \Cc'$  by
 $\langle T_2^*\omega,\ell\rangle=\langle\omega_A,\ell_A\rangle+\langle\omega_B,\ell_B\rangle
 :=\langle\omega_A\oplus\omega_B,\ell\rangle.$ That is
\begin{equation}\label{eq-120}
    T_2^*\omega=\omega_A\oplus\omega_B\in \Cc^{\prime *}.
\end{equation}
  This yields
 \begin{equation*}
    \Lb(\omega)=\delta_{\overline{\Gamma}}(\omega_A\oplus\omega_B),\quad
    \omega\in\Yii
\end{equation*}
and the extended dual problem \Dii\ is
\begin{equation*}
   \textsl{maximize }
   \langle\omega_A,\mu\rangle+\langle\omega_B,\nu\rangle ,\quad
   \omega\in\Yii \textsl{ such that } \omega_A\oplus\omega_B\in\overline{\Gamma}
   \tag{$D_2$}
\end{equation*}
Note that for this dual problem to be meaningful, it is necessary
that \HC\ holds: i.e. $(\mu,\nu)\in\Xi.$ This is realized if
$(\mu,\nu)\in\dom\Ls$ or equivalently if $\inf(\MK)<\infty.$

\subsubsection*{The constraint qualification}
One will be allowed to apply Theorem \ref{T3b} under the
constraint qualification (\ref{xCQbis}):
\begin{equation}\label{eq-20}
    (\mu,\nu)\in \diffdom \Ls.
\end{equation}
Let us give some details on this abstract requirement.

\begin{remark}\label{rem-04}
Note that for all $\mu\in\PA,$ $\nu\in\PB,$
$(\mu,\nu)\not\in\icordom\Ls$ if $\AB$ is an infinite set. Indeed,
for all $\pi\in P(\mu,\nu)$ such that $\IAB c\,d\pi<\infty,$ one
can find $(a_o,b_o)$ such that with $\varepsilon_{(a_o,b_o)}$ the
Dirac measure  at $(a_o,b_o),$
$\ell_t:=t\varepsilon_{(a_o,b_o)}+(1-t)\pi\not\geq 0$ for all
$t<0,$ so that $\Fs_1(\ell_t)=+\infty$ (Proposition
\ref{res-02}-a). This shows that
$[\ell_0,\ell_1]=[\pi,\varepsilon_{(a_o,b_o)}]\subset \dom\Fs_1$
while $\ell_t\not\in \dom\Fs_1$ for all $t<0.$ Hence,
$(\mu,\nu)\not\in\icordom\Ls$ and one has to consider the
assumption (\ref{eq-20}) on $(\mu,\nu)$ rather than
$(\mu,\nu)\in\icordom\Ls.$
\\
This is in contrast with the situation encountered in
\cite{Leo06b} where the rule is $x_o\in\icordom\Ls.$
\end{remark}

\begin{lemma}\label{res-03}We have $\dom\Ls=\diffdom\Ls.$
\end{lemma}
 \proof  Proposition \ref{res-02}-a
states that $\dom\Fs_1\subset \Lp$ where $\Lp=\{\ell\in \Cc';
\ell\geq 0\}$ is the the \positif\  cone of $\Cc'.$ Therefore,
$\Fs_1=\Fs_1+\delta_{\Lp}.$ Consequently, with (\ref{xped}) one
obtains that
\[
\Ls(x)=\inf\{\Fs_1(\ell); \ell\in \Lp, T\ell=x\},\quad x\in\Xi.
\]
Suppose ad absurdum that there is some $x_o\in \dom\Ls$ such that
$x_o\not\in\diffdom\Ls.$ This implies that there exists some
half-line $]x_o,x_o+\infty(x_o'-x_o)[$ on which $\Ls$ achieves the
value $+\infty,$ which in turn implies that $\Fs_1$ must achieve
the value $+\infty$ somewhere on $\Lp.$ But this is impossible
since $\Fs_1(\ell)=\|\ell\|_c^*$  for all $\ell\in \Lp.$ This
completes the proof of the lemma.
\endproof

As a consequence of this lemma, it appears that (\ref{eq-20}) is
\emph{not} a constraint qualification. One can apply Theorem
\ref{T3b} under the only restriction that $\inf(\MK)<\infty.$ This
gives the following

\begin{lemma}\label{res-17}
Let us assume that $\inf(\MK)<\infty.$ Then, $(P)$ and \Dii\ both
admit a solution in $\PAB$ and $\Yii.$ Furthermore, any
$(\pi,\omega)\in\PAB\times \Yii$ is a solution of $(P)$ and \Dii\
if and only if
\begin{equation}\label{eq-17}
\left\{\begin{array}{ll}
  (a) & \pi_A=\mu, \pi_B=\nu, \IAB c\,d\pi<\infty;\\
  (b) & \pi\in \partial_{\Cc'}\Fd(\eta) \textsl{ where } \\
  (c) & \eta=T_2^\ast \omega. \\
\end{array}\right.
\end{equation}
\end{lemma}
As $\Fs_1$ and $\Fd$ are mutually convex conjugates,
(\ref{eq-17}-b) is equivalent to
\begin{equation}\label{eq-122}
    \eta\in \partial_{\Cc^{\prime *}}\Fs_1(\pi)
\end{equation}
and also equivalent to Young's identity
\begin{equation}\label{eq-122b}
    \Fs_1(\pi)+\Fd(\eta)=\langle\eta,\pi\rangle
\end{equation}
and also equivalent to
\begin{equation*}
    \left\{\begin{array}{l}
      \Fd(\eta)= 0\\
      \langle\eta,\pi\rangle=\IAB cd\pi.\\
    \end{array}\right.
\end{equation*}
In other words:
\begin{theorem}\label{res-20}
Let $\pi\in P(\mu,\nu)$ be such that $\IAB c\,d\pi<\infty.$ Then:
\begin{enumerate}
\item \Dii\ admits at least a solution in $\Yii;$

\item $\pi$  is an optimal plan if and only if there exists some
$\omega\in\Yii$ such that
\begin{equation*}
    \left\{\begin{array}{cl}
         (a) & T_2^*\omega\in\overline{\Gamma}  \\
     (b) & \langle\omega,(\mu,\nu)\rangle=\IAB cd\pi;\\
    \end{array}\right.
\end{equation*}
With $\eta=T_2^*\omega,$ this implies the equivalent statements
(\ref{eq-17}-b), (\ref{eq-122}) or (\ref{eq-122b}).

\item If such an $\omega$ exists, it is a solution of \Dii\ and
any other solution of \Dii\ is also convenient.
\end{enumerate}
\end{theorem}
This is the core of the extended saddle-point method applied to
Monge-Kantorovich problem. To prove a practical optimality
criterion one still has to translate these abstract properties.

\subsection{Optimal plan: preliminary results for the sufficient
condition}\label{sec-02}

The next lemmas are preliminary results for the proof of a
sufficient condition for the optimality.

\begin{lemma}\label{res-18}
Let $\varphi$ and $\psi$ be real functions on $A$ and $B.$
\begin{enumerate}
    \item The \lsc\ regularizations $\ls\varphi$ and $\ls\psi$ of $\varphi$ and $\psi$ satisfy
    $$
    \ls(\varphi\oplus\psi)=\ls\varphi\oplus\ls\psi.
    $$
    \item If $\varphi$ and $\psi$ are such
that $\varphi\oplus\psi=c$ on some subset $\mathcal{T}$ of $\AB$
and $\varphi\oplus\psi\leq c$ everywhere on $\AB.$ Then,
$\ls\varphi$ and $\ls\psi$ still share the same properties.
\end{enumerate}
\end{lemma}

\proof \boulette{1} For each $(a,b)\in\AB,$
\begin{eqnarray*}
   \ls (\varphi\oplus\psi)(a,b)
  &=& \sup_{V\in\mathcal{N}((a,b))}\inf_{(a',b')\in V}[\varphi(a')+\psi(b')] \\
  &=& \sup\left\{\inf_{(a',b')\in V_A\times V_B}[\varphi(a')+\psi(b')];
  V_A\in\mathcal{N}(a),V_B\in\mathcal{N}(b)\right\} \\
  &=& \sup\left\{\inf_{a'\in V_A} \varphi(a'); V_A\in\mathcal{N}(a)\right\}
      + \sup\left\{\inf_{b'\in V_B} \psi(b'); V_B\in\mathcal{N}(b)\right\}\\
  &=& \ls\varphi(a)+\ls\psi(b)
\end{eqnarray*}
where $\mathcal{N}(x)$ stands for the set of all open
neighbourhoods of $x.$

\Boulette{2} It is a direct consequence of the lower
semicontinuity of $c$ and statement (1).
\endproof

\begin{lemma}\label{res-19}
Let $\pi\in P(\mu,\nu)$ be such that $\IAB c\,d\pi<\infty.$
Suppose that there exists two real-valued functions
$\varphi\in\mathbb{R}^A$ and $\psi\in\mathbb{R}^B$ such that
\begin{equation}\label{eq-82}
    \left\{\begin{array}{ll}
      \varphi\oplus\psi\leq c &\ \textrm{everywhere} \\
       \varphi\oplus\psi= c &\ \pi\textrm{-almost everywhere}.\\
    \end{array}\right.
\end{equation}
\begin{enumerate}
    \item If $\varphi$ is $\mu$-measurable and $\psi$ is
    $\nu$-measurable, then $\varphi\in L_1(A,\mu)$ and $\psi\in L_1(B,\nu).$
    \item In any case, the real-valued functions $\ls\varphi$ and
$\ls\psi$ still satisfy (\ref{eq-82}) together with $\ls\varphi\in
L_1(A,\mu)$ and $\ls\psi\in L_1(B,\nu).$
\end{enumerate}
\end{lemma}

\proof \boulette{1} Let us fix $(a_o,b_o)\in\AB$ such that
$c(a_o,b_o)<\infty$ (such a point exists since $\IAB
c\,d\pi<\infty$ for some $\pi.$) We have
$\varphi(a)=c(a,b_o)-\psi(b_o)\geq -\psi(b_o)$ for all $a\in A$
and similarly $\psi\geq -\varphi(a_o).$  Hence, the integrals
$\int_A \varphi\,d\mu\in [-\psi(b_o),+\infty]$ and
$\int_B\psi\,d\nu\in[-\varphi(a_o),+\infty]$ are well-defined.
Finally, $\varphi\in L_1(A,\mu)$ and $\psi\in L_1(B,\nu)$ since
$\int_A \varphi\,d\mu+\int_B\psi\,d\nu= \IAB c\,d\pi<\infty.$

\Boulette{2} Applying Lemma \ref{res-18} with $\mathcal{T}$ a
measurable set such that $\pi(\mathcal{T})=1$ yields two lower
bounded measurable functions $\ls\varphi$ and $\ls\psi$ which
still satisfy (\ref{eq-82}). One concludes as above.
\endproof

Let $\overline{\Upsilon}$ be the $\sigma(\Yii,\Xi)$-closure of
\begin{equation}\label{eq-83}
    \Upsilon=\{(\varphi,\psi)\in C_A\times C_B;
\varphi\oplus\psi\leq c\}.
\end{equation}
\begin{lemma}\label{res-14}\
\begin{enumerate}[(a)]
    \item For all $(a,b)\in\s,$ $\Ls(\varepsilon_a,\varepsilon_b)=c(a,b).$
    \item For any  $\omega\in\Yii,$ we have $\omega\in \overline{\Upsilon}$
if and only if $\langle\omega,\kappa\rangle\leq \Ls(\kappa),
\forall \kappa\in\Xi.$
    \item $T_2^*\overline{\Upsilon}\subset\overline{\Gamma}.$
\end{enumerate}
\end{lemma}

\proof \boulette{a} For any $(a,b)\in\s,$
$\Ls(\varepsilon_a,\varepsilon_b)=\inf\{\IAB
c\,d\pi;\pi\in\PAB:\pi_A=\varepsilon_a,\pi_B=\varepsilon_b\}=\IAB
c\, d\varepsilon_{(a,b)}=c(a,b)$ where we used the dual equality
(\ref{xped}) and the fact that $\varepsilon_{(a,b)}$ is the unique
plan $\pi$  with marginals $\varepsilon_a$ and $\varepsilon_b.$

\Boulette{b} It is enough to check that for all
$\phi=(\varphi,\psi)$ in $C_A\times C_B$
\begin{equation}\label{eq-126}
    \phi\in \Upsilon\Leftrightarrow [\langle\phi,\kappa\rangle\leq
    \Ls(\kappa), \forall \kappa\in\Xi].
\end{equation}
 Young's inequality $\langle\phi,\kappa\rangle\leq
\La(\phi)+\Ls(\kappa), \forall \phi,\kappa$ and
$\phi\in\Upsilon\Leftrightarrow\La(\phi)=\delta_\Upsilon(\phi)=0$
 give the direct implication. For the converse, choosing
$\kappa=(\varepsilon_a,\varepsilon_b)$  in the right-hand side of
(\ref{eq-126}), one obtains with the previous statement (a) that
$\varphi\oplus\psi\leq c.$

\Boulette{c} It is clear that $T^*\Upsilon\subset\Gamma$ and one
concludes with the $\sYXii$-$\sULii$-continuity of $T_2^\ast:
\Yii\to\Uii,$ see Lemma \ref{L2}-d.
\endproof

\begin{lemma}\label{res-13}
Let $\pi\in P(\mu,\nu)$ be such that $\IAB c\,d\pi<\infty$ and
suppose that there exist two real functions $\varphi$ in
$L_1(A,\mu)$ and $\psi$ in $L_1(B,\nu)$ satisfying (\ref{eq-82}).
\\
Then, there exists some $\omega$ in $\overline{\Upsilon}$ such
that
\begin{equation*}
    \left\{\begin{array}{l}
      \langle\omega,(\mu,\nu)\rangle=\IAB c\,d\pi; \\
      \omega(\epsilon_a,\epsilon_b)=\varphi(a)+\psi(b)\textrm{ for $\pi$-a.e. $(a,b)$ and} \\
      \omega(\kappa)\leq \Ls(|\kappa|),\ \forall \kappa\in \Xi.\\
    \end{array}\right.
\end{equation*}
\end{lemma}

\proof There exists a measurable subset $\mathcal{T}$ of $\s$ such
that $\pi(\mathcal{T})=1$ and $\varphi\oplus\psi=c$ everywhere on
$\mathcal{T}.$ Let $E_o$ be the vector subspace of $\Xi$ spanned
by $(\mu,\nu)$ and $\{(\epsilon_a,\epsilon_b); (a,b)\in
\mathcal{T}\}.$ It follows from our assumptions on $\varphi$ and
$\psi$ that for all \positif\ $\kappa=(\kappa_1,\kappa_2)\in E_o,$
${\varphi}$ is in $L_1(A,\kappa_1)$ and ${\psi}$ is in
$L_1(B,\kappa_2).$ Define the linear form $\omega_o$ on $E_o$ for
each $\kappa\in E_o$ by
$$
\omega_o(\kappa)=\int_A\varphi\,d\kappa_1 +\int_B\psi\,d\kappa_2.
$$
Clearly,
\begin{equation}\label{eq-81}
    \omega_o(\mu,\nu)=\IAB c\,d\pi
\end{equation}
and  for all \positif\ $\kappa\in E_o,$
\begin{eqnarray*}
  \omega_o(\kappa)
  &=& \int_A\varphi\,d\kappa_1+\int_B\psi\,d\kappa_2\\
   &\leq& \sup\left\{\int_A\tilde{\varphi}\,d\kappa_1+\int_B\tilde{\psi}\,d\kappa_2;\
    \tilde{\varphi}\in  L_1(A,\kappa_1),\ \tilde{\psi}\in
    L_1(B,\kappa_1),\
   \tilde{\varphi}\oplus\tilde{\psi}\leq c\right\}.\\
\end{eqnarray*}
Denoting $(K_{\kappa})$ and $(\overline{K}_{\kappa})$ the
analogues of problems $(K)$ and $(\overline{K})$ with
$(\kappa_1,\kappa_2)$ instead of $(\mu,\nu),$ this means that
$$
\omega_o(\kappa)\leq\sup(\overline{K}_{\kappa}).
$$
The dual equality (\ref{xed1}) states that
$\sup(K_{\kappa})=\Ls({\kappa}).$ As we have already seen at
Theorem \ref{res-MK1}-a that
$\sup(\overline{K}_{\kappa})=\sup(K_{\kappa}),$ we obtain:
$\sup(\overline{K}_{\kappa})=\Ls({\kappa}).$ Therefore, we have
proved that $\omega_o(\kappa)\leq\Ls(\kappa),$ for all $\kappa\in
E_o,$ $\kappa\geq 0.$ As for any $\kappa\in E_o,$
$\omega_o(\kappa)=\IAB c\,d\rho$ for any measure $\rho$ with
marginals $\rho_A=\kappa_1$ and $\rho_B=\kappa_2,$ one sees that
$\omega_o$ is \positif. It follows that
$$
\omega_o(\kappa)\leq\Ls(|\kappa|),\quad\kappa\in E_o
$$
where $|\kappa|=(|\kappa_1|,|\kappa_2|)$ and $|\kappa_i|$ is the
absolute value of the measure $\kappa_i.$
\\
Note that $\Xi$ is a Riesz space since it is the topological dual
of a normed Riesz space. Hence, any $\kappa\in\Xi$ admits
\positif\ and \negatif\ parts $\kappa^+$ and $\kappa^-$, and its
absolute value is $|\kappa|=\kappa^++\kappa^-.$ This allows to
consider the positively homogeneous convex function
$\Ls(|\kappa|)$ on the vector space $E_1$ spanned by $\dom\Ls.$ By
the analytic form of Hahn-Banach theorem, there exists an
extension $\omega$ of $\omega_o$ to  $E_1$ which satisfies
$\omega(\kappa)\leq \Ls(|\kappa|)$ for all $\kappa\in E_1.$ But
$E_1=\Xi$ and one completes the proof of the lemma with
(\ref{eq-81}) and Lemma \ref{res-14}-b.
\endproof

\subsection{Optimal plan: preliminary results for the necessary
condition}\label{sec-03}

Under the condition (\ref{eq-17}-a), $\pi$ necessarily satisfies:
$\suppp\subset\s.$ This fact will be invoked without warning.

\begin{lemma}\label{res-10a}
Let $\pi$ and $\eta$ satisfy (\ref{eq-17}-a,b). Then, the
restriction of $\eta$ to $L_\infty.\pi:=\{h.\pi; h\in
L_\infty(A\times B,\pi)\}$ is given by
\begin{equation}\label{eq-121}
    \langle \eta,h.\pi\rangle
    =\IAB h c\,d\pi, \quad \forall h\in L_{\infty}(\pi)
\end{equation}
\end{lemma}

 \proof
To specify the restriction $\gamma$ of $\eta$ to $L_{\infty}.\pi,$
it is enough to vary $\Fs_1$ in the direction $L_{\infty}.\pi$ to
get with (\ref{eq-122}):
$\gamma\in\partial_{(L_{\infty}.\pi)^*}\Fs_1(\pi).$ Taking $h\in
L_{\infty}(\pi)$ such that $\|h\|_\infty\leq 1,$ by monotone
convergence we obtain $\Fs_1(\pi+h.\pi)=\sup\{\IAB (1+h)u\,d\pi;
u\in \Cc, u\leq c\}=\IS (1+h)c\,d\pi.$ It comes out that
$\partial_{(L_{\infty}.\pi)^*}\Fs_1(\pi)=\{c\},$ which gives
(\ref{eq-121}).
\endproof

We first derive the necessary condition in the special case where
$c$ is assumed to be finite and continuous.

\begin{proposition}\label{res-11b}
Assume that $c$ is finite and continuous and let $\pi$ be an
optimal plan. Then, there exist  two finitely-valued \usc\
functions $\varphi$ on $A$ and $\psi$ on $B$ such that
\begin{equation*}
    \left\{\begin{array}{ll}
      \varphi\oplus\psi\leq c & \textrm{everywhere and} \\
       \varphi\oplus\psi= c &\textrm{on } \suppp. \\
    \end{array}\right.
\end{equation*}
\end{proposition}
\proof At the beginning of this proof, $c$ is only assumed to be
finite and \lsc. By Lemma \ref{res-03}, (\ref{eq-20}) is
satisfied. Let $\eta$ and $\omega$ be as in (\ref{eq-17}-b \& c).
Because of Theorem \ref{T3b}-c-1 \& 2, there exists a generalized
sequence $\{(\alpha_\tau,\beta_\tau)\}$ in $\dom \La_1$ such that
$\lim_\tau T_1^*(\alpha_\tau,\beta_\tau)=\eta$ with respect to
$\sigma(\Cc^{\prime *},\Cc').$ As $T_1^*\dom\La_1\subset \Ui$ (see
Lemma \ref{L2}-g),  $T_1^*(\alpha_\tau,\beta_\tau)=
\alpha_\tau\oplus\beta_\tau\in \Cc$ and
\begin{equation}\label{eq-64}
\left\{\begin{array}{ll}
  (a) & \lim_\tau \alpha_\tau\oplus\beta_\tau=\eta \textrm{ with } \\
  (b) & \Cc\ni\alpha_\tau\oplus\beta_\tau\leq c\textrm{ for all }\tau \\
\end{array}\right.
\end{equation}
Defining
\begin{equation}\label{eq-63}
    \etat(a,b)=\langle\eta,\epsilon_{(a,b)}\rangle,\quad
    (a,b)\in\AB,
\end{equation}
where $\epsilon_{(a,b)}$ is the Dirac mass at $(a,b),$ one
immediately sees that
\begin{equation}\label{eq-61}
    \etat\leq c.
\end{equation}
Furthermore, since $\etat=\lim_\tau \alpha_\tau\oplus\beta_\tau$
pointwise ($\Cc'$ contains the Dirac masses), by (\cite{BL92},
Corollary 3.5) we obtain that
$
 \etat=\varphi\oplus\psi
$
for some functions $\varphi$ and $\psi$ on $A$ and $B.$ This gives
us some hope to complete the proof, but as will be seen below,
$\etat$ isn't the right function to be considered.
\\
For any $(a,b)$ in  $\supp\pi:$ the support of $\pi,$ one can find
a sequence $\{h_k\}$ in $\Cc$ such that $\lim_k
h_k.\pi=\epsilon_{(a,b)}$ in $\Pc,$ see Lemma \ref{res-12} below.
As $c$ is \lsc, with (\ref{eq-121}) we obtain
\begin{equation}\label{eq-60}
    \liminf_k \langle h_k.\pi,\eta\rangle
 =\liminf_k\IAB h_k c\,d\pi \geq c(a,b).
\end{equation}
 Unfortunately,  no
 regularity property for $\eta$ has been established to insure
 that $\langle \eta,\epsilon_{(a,b)}\rangle\geq \liminf_k
\langle\eta, h_k.\pi\rangle;$ this would lead to the converse of
(\ref{eq-61}): $\etat\geq c$ on $\suppp.$ An alternate strategy is
to introduce the \usc\ regularization
 \begin{equation*}
    \etab=\us\etat
\end{equation*}
of $\etat$ on $\AB.$ As $\etab$ is \usc, for all
$(a,b)\in\supp\pi,$ we have $\etab(a,b)\geq\limsup_k\IAB \etab
h_k\,d\pi.$ Now, one obtains with (\ref{eq-60}) that
\begin{equation}\label{eq-65a}
      \etab(a,b)\geq c(a,b),\quad \forall (a,b)\in\suppp.
\end{equation}
Regularizing both sides of (\ref{eq-61}) and assuming that $c$ is
\usc\ and therefore \emph{continuous}, we obtain that
\begin{equation}\label{eq-65b}
    \etab\leq c.
\end{equation}
It remains to check that
\begin{equation*}
    \etab=\varphi\oplus\psi
\end{equation*}
for some finitely-valued \usc\ functions $\varphi$ and $\psi$ on
$A$ and $B.$ With (\ref{eq-120}) and (\ref{eq-17}) we know that
$\eta=\omega_A\oplus\omega_B$ for some $\omega\in\Yii.$ It follows
that $\etat=\widetilde{\omega}_A\oplus\widetilde{\omega}_B$ where
$\widetilde{\omega}_A(a)=\omega_A(\epsilon_a)$ and
$\widetilde{\omega}_B(b)=\omega_B(\epsilon_b).$ With Lemma
\ref{res-18}, one sees that
$\etab=\us\widetilde{\omega}_A\oplus\us\widetilde{\omega}_B.$ This
proves the desired result with $\varphi=\us\widetilde{\omega}_A$
and $\psi=\us\widetilde{\omega}_B.$ Since $\etat\leq\etab\leq c$
and both $\etat$ and $c$ are finitely-valued, so are $\varphi$ and
$\psi.$
\endproof

\begin{remark}
    By means of the usual approaches \cite{Rus96,AP02,Vill03}, one
    can prove when $c$ is finitely-valued that under the assumptions (\ref{eq-29a}) or (\ref{eq-29b}),
    $\varphi$ and $\psi$ can be required to be $c$-concave conjugates
    to each other.
    In the special case where $c$ is assumed to be continuous,
    $c$-concave conjugates
    are \usc. This is in accordance with Proposition \ref{res-11b}.
\end{remark}

Now, we remove the assumption that $c$ is finite and continuous
and only assume that it is \lsc. The main technical result for the
proof of the characterization of the optimal plans is the
following

\begin{lemma}\label{res-11c}
Assume that $c$ is a $[1,\infty]$-valued \lsc\ function.  Let
$\pi\in\Pc$ and $\eta\in\Cc^{\prime *}$ be as in (\ref{eq-17}-b),
i.e. $
     \pi\in \partial_{\Cc'}\Fd(\eta)
$ and  define the function $\etat$ on $\s$ by
\begin{equation}\label{eq-71}
    \etat(a,b)=\langle\eta,\epsilon_{(a,b)}\rangle,\quad
    (a,b)\in\s.
\end{equation}
Then,
\begin{equation*}
    \left\{\begin{array}{ll}
     \etat\leq c&\ \textrm{on } \s\\
     \etat=c&\ \textrm{on }\supp\pi\bigcap\{c<\infty\}\\
    \end{array}\right.
\end{equation*}
and $\etat$ is a finitely-valued measurable function on $\s.$
\end{lemma}
\begin{remark}
We assume that $c\geq 1$ without loss of generality, see Remark
\ref{rem-03}, to allow dividing by $c$ in the definition of $\Cc.$
\end{remark}
 \proof
Because of Theorem \ref{T3b}-c-2, there exists a sequence
$\{\rho_n\}$ in $\CAB$ such that
\begin{equation}\label{eq-49}
\left\{\begin{array}{l}
  \rho_n\leq c, \forall n \textrm{\quad and} \\
   \lim_{n\rightarrow\infty}\rho_n=\eta \\
\end{array}\right.
\end{equation}
with respect to $\sigma(\Cc^{\prime*},\Cc').$ Having Remark
\ref{rem-02} in mind, recall that only the restriction of $\rho_n$
to $\s$ carries information as regards to the dual pairing
$\langle\Cc^{\prime*},\Cc'\rangle.$ Also recall that the
$k^{\mathrm{th}}$\emph{ Moreau-Yosida approximation} of a function
$u$ on a space with metric $d$ is defined for all $x$ by
$u^{(k)}(x)=\inf_y\{u(y)+kd(x,y)\}.$ Defining the Moreau-Yosida
approximations
\begin{eqnarray*}
      \rho_{n,k}&=&[\max(\rho_n,k)]^{(k)} \\
      c_k&=&[\max(c,k)]^{(k)}\\
\end{eqnarray*}
for all $n,k\geq1,$ (\ref{eq-49}) implies that
\begin{equation}\label{eq-47}
    \left\{\begin{array}{cl}
      (a)&\rho_{n,k}, c_k\in \CAB, \ \rho_{n,k}\leq c_k\leq c, \ \forall n\geq 1 \\
      (b)&\lim_n\lim_k\langle\rho_{n,k},m\rangle=\langle\eta,m\rangle,\ \forall m\in\M, m\geq 0\\
      (c)&0\leq c_k\uparrow c \textrm{ pointwise.}\\
    \end{array}\right.
\end{equation}
where $\M$ is the space of all measures $m$ on $\s$ such that $\IS
c\,d|m|<\infty.$ By Remark \ref{rem-02}, one sees that
$\M\subset\Cc'.$
\\
 While deriving (\ref{eq-47}), we used the well-known results:
\begin{itemize}
    \item a Moreau-Yosida approximation is a continuous function and
    \item the sequence of
Moreau-Yosida approximations of a function tends pointwise and
\croissant ly to its lower semicontinuous regularization.
\end{itemize}
The proof of statement (\ref{eq-47}-b) relies on the monotone
convergence theorem; this is the reason why it holds for all $m$
in $\M$ rather than in $\Cc'.$

Let us introduce the cone $\Qp=\{\ell\in \Cc'; \ell\geq 0 \textrm{
and } \langle \eta,\ell\rangle\geq 0\}$ and $\Q$ the vector space
spanned by $\Qp.$ We first consider the restriction $\theta$ of
$\eta$ to $\Q.$ By (\ref{eq-121}), $\pi$ is in $\Qp$ and
(\ref{eq-122}) gives us $\theta\in\partial_{\Q^*}N(\pi)$ where $
N(\ell)=\sup\{\langle u,\ell\rangle; u\in \Cc, |u|\leq c\},$
$\ell\in \Q $ which is the dual norm $\|\cdot\|_c^*$ restricted to
$\Q.$ It follows that $\theta$ belongs to the topological dual
space $\Q'$ of the normed space $(\Q,\|\cdot\|_c^*):$
\begin{equation*}
    \theta:=\eta_{|\Q}\in\Q'.
\end{equation*}
This topological regularity of $\theta$  will allow us a few lines
below to invoke Br\o nsted-Rockafellar's lemma. It is not clear
that $\eta$ is continuous on the whole normed space $\Cc'.$
\\
Let us denote $\Psi$ the restriction of $\Fs_1$ to $\Q$ and
$\Psi^*$ its convex conjugate with respect to the dual pairing
$\langle\Q,\Q'\rangle.$  Since $\theta\in \Q',$ $\theta\geq 0$ and
$\eta\in\overline{\Gamma},$ one sees that
$0\leq\Psi^*(\theta)\leq\Fd(\eta)=0.$   As (\ref{eq-17}-b) is
equivalent to Young's identity (\ref{eq-122b}), one obtains
\begin{equation*}
    \Psi(\pi)+\Psi^*(\theta)=\langle\theta,\pi\rangle=\lim_{n\rightarrow\infty}\langle\xi_n,\pi\rangle
\end{equation*}
where
\begin{equation}\label{eq-46}
    \xi_n={\rho_{n,k(n)}}_{|\Q}\in\Q'
\end{equation}
is the restriction of $\rho_{n,k(n)}\in \Cc^{\prime *}$ to
$\Q\subset \Cc'$ for some sequence $\{k(n)\}_n$ which converges
fast enough to infinity to imply that
$\lim_{n\rightarrow\infty}\langle\xi_n,\pi\rangle=\langle\theta,\pi\rangle$
by means of (\ref{eq-47}-b).

Denote $\Psi_n$  the restriction to $\Q$ of the analogue of
$\Fs_1$ with $c_{k(n)}$ instead of $c$ and $\Psi^*_n$  its convex
conjugate with respect to $\langle\Q,\Q'\rangle.$ By
(\ref{eq-47}-c), we have
$\lim_{n\rightarrow\infty}\Psi_{n}(\pi)=\Psi(\pi).$  By
(\ref{eq-47}-a), we also have $\Psi^*(\theta)=\Psi^*_{n}(\xi_n)=0$
for all $n.$  Therefore,
\begin{equation*}
    \Psi_n(\pi)+\Psi_n^*(\xi_n)=\langle\xi_n,\pi\rangle+\epsilon_n
\end{equation*}
with $\lim_{n\rightarrow\infty}\epsilon_n=0.$ In other words,
$\xi_n$ is an $\epsilon_n$-subgradient of $\Psi_n$ at $\pi.$
Hence, by Br\o nsted-Rockafellar lemma, there exist two sequences
$\{\pi_n\}$ in $\Q$ and $\{\theta_n\}$ in $\Q'$ such that for all
$n,$
\begin{align}
    &\|\pi_n-\pi\|\leq \sqrt{\epsilon_n},\label{eq-45a}\\
    &\|\xi_n-\theta_n\|\leq\sqrt{\epsilon_n}\label{eq-45b}
\end{align}
(both norms $N(\ell)=\sup\{\langle u,\ell\rangle; u\in \Cc,
|u|\leq c\}$ on $\Q$ and
$\sup\{\langle\cdot,\ell\rangle;\ell\in\Q,N(\ell)\leq 1\}$ on
 $\Q'$ are simply written $\|\cdot\|$) and
\begin{equation}\label{eq-50}
   \theta_n\in\partial_{\Q'}\Psi_n(\pi_n).
\end{equation}
We define
\begin{equation*}
    G=\{(a,b)\in\s;\langle\eta,\epsilon_{(a,b)}\rangle\geq 0\}
\end{equation*}
the set of all $(a,b)\in\AB$ such that $\epsilon_{(a,b)}\in\Q.$
Since $c_{k(n)}$ is finite and continuous, proceeding as in
Proposition \ref{res-11b}, one shows as for (\ref{eq-65a}) that
\begin{equation}\label{eq-43}
   \thetab_n(a,b)=c_{k(n)}(a,b),\quad \forall (a,b)\in\cl G\cap\supp\pi_n
\end{equation}
where
\begin{equation*}
    \thetab_n(a,b)=\us\thetat_n(a,b),\quad (a,b)\in \s
\end{equation*}
is the \usc\ regularization of
$$
\thetat_n(a,b)=\left\{\begin{array}{ll}
  \theta_n(\epsilon_{(a,b)}) & \textrm{if }(a,b)\in G \\
  -\infty & \textrm{otherwise } \\
\end{array}\right.,\quad
   (a,b) \in \s
 $$
and $\cl G$ is the closure of $G$ in $\AB.$ Since $\s$ is closed,
we have $\cl G\subset \s.$
 As $\pi_n$ may not be a measure, one uses Lemma \ref{res-12}
below instead of its usual analogue.

Thanks to (\ref{eq-45a}), $\limn\pi_n=\pi$ \emph{strongly} in $\Q$
 and for all large enough $n$ we have $ \pi_n\in\Qp$ and
\begin{equation}\label{eq-67a}
      \supp\pi_n=\suppp.
\end{equation}
With $\etab$  the \usc\ regularization of $\etat,$ we have
$\{\etab\geq 0\}=\cl\{\etat\geq 0\}.$  It follows from
(\ref{eq-65a}) that
\begin{equation}\label{eq-69}
    \suppp\subset \cl G.
\end{equation}
 Thanks to (\ref{eq-46}) and (\ref{eq-45b}), for all  $0\leq r<\infty,$
\begin{equation}\label{eq-67b}
    \limn \sup_{(a,b)\in G\cap \{c\leq
    r\}}|\rho_{n,k(n)}(a,b)-\thetat_n(a,b)|=0.
\end{equation}
Let us assume for a while that
\begin{equation}\label{eq-68}
    \sup_{\s}c<\infty.
\end{equation}
Under this assumption, (\ref{eq-67b}) leads us to  $ \limn
\sup_G|\rho_{n,k(n)}-\thetat_n|=0$ on $G.$  Upper regularizing,
because of this uniform estimate and the continuity of
$\rho_{n,k(n)},$ one obtains
\begin{equation}\label{eq-70}
    \limn \gamma_n=0
\end{equation}
where $ \gamma_n=\sup_{\cl G}|\rho_{n,k(n)}-\thetab_n|.$ By
(\ref{eq-43}), (\ref{eq-67a}), (\ref{eq-69}) and $\rho_{n,k}\leq
\rho_n$ for all $n,k,$ we obtain for all large enough $n:$
$\rho_{n}(a,b)\geq c_{k(n)}(a,b)-\gamma_n,$ $\forall
(a,b)\in\suppp.$ Letting $n$ tend to infinity, we see with
(\ref{eq-49}), (\ref{eq-47}-c) and (\ref{eq-70}) that
$\etat(a,b)=c(a,b),$ for all $(a,b)\in\suppp$ where $\etat$ is
defined at (\ref{eq-71}). We have just proved that under the
assumption (\ref{eq-68}),
\begin{equation}\label{eq-73}
    \left\{\begin{array}{ll}
     (a)\quad \etat(a,b)\leq c(a,b),&\ \forall (a,b)\in\s\\
     (b)\quad \etat(a,b)=c(a,b),&\ \forall (a,b)\in\suppp\\
    \end{array}\right.
\end{equation}
where the statement (a) directly follows from  (\ref{eq-49}).

It remains to remove the restriction (\ref{eq-68}). For each
$k\geq 1,$ let
\begin{equation*}
    \left\{\begin{array}{lll}
      \s_k & = & \{c\leq k\} \textrm{\quad and} \\
      c_k & = & c+\delta_{\s_k} \\
    \end{array}\right.
\end{equation*}
The function $c_k$ is \lsc\  on $\AB$ and satisfies (\ref{eq-68});
 $\{\s_k\}$ is a \croissant\ sequence of  closed level sets of $c$
  with $\s_k\subset\s$ for all $k.$  By Proposition \ref{res-02}-b we have
$\suppp\subset\s.$
\\
It is assumed that $\pi\in
\partial_{\Cc'}\Fd(\eta)$ which is equivalent to the Young's identity
$\Fs_1(\pi)+\Fd(\eta)=\langle\eta,\pi\rangle$ or equivalently $[
  \Fs_1(\pi) = \langle\eta,\pi\rangle \textrm{ and }
  \Fd(\eta) = 0]$ which is also equivalent to
\begin{equation}\label{eq-74}
\left\{\begin{array}{l}
  \langle\eta,\pi\rangle =  \IAB c\,d\pi \quad $and$ \\
   \langle\eta,\ell\rangle\leq \Fs_1(\ell), \ \forall\ell\in\Cc'\\
\end{array}\right.
\end{equation}
because of  Proposition \ref{res-02}-d and Lemma \ref{res-14a}
below. Let us consider for each $k$
\begin{equation*}
\left\{\begin{array}{rcl}
  \pi_k & = & \1_{\s_k}.\pi \\
  \langle\eta_k,\ell\rangle & = & \langle\eta,\ell\rangle,
  \quad \forall \ell\in\Cc' \textrm{\ such that\ }\supp\ell\subset\s_k \\
\end{array}\right.
\end{equation*}
Note with Remark \ref{rem-02} that $\pi_k\in C_{c_k}'$ and
$\eta_k\in C_{c_k}^{\prime *}.$ Also introduce $\Theta_k$ and
$\Theta^*_k$ the analogues of $\Fs_1$ and $\Fd$ where $c$ is
replaced by $c_k.$
\\
By Proposition \ref{res-02}-d and (\ref{eq-121}),
$\Theta_k(\pi_k)=\IAB c_k\,d\pi_k=\langle\eta_k,\pi_k\rangle.$
Since for any $\ell\in\Cc'$ such that $\supp\ell\subset\s_k $ we
have $\Theta_k(\ell)=\Fs_1(\ell),$ one obtains with (\ref{eq-74})
that $\langle\eta_k,\ell\rangle\leq\Theta_k(\ell)$ for all
$\ell\in C_{c_k}'.$ Reasoning as for the derivation of
(\ref{eq-74}) but taking the reverse way, this shows that
$$\pi_k\in\partial_{C_{c_k}'}\Theta^*_k(\eta_k).$$
Applying (\ref{eq-73}) yields
\begin{equation*}
    \left\{\begin{array}{ll}
     \etat_k(a,b)\leq c_k(a,b),&\ \forall (a,b)\in\s_k\\
     \etat_k(a,b)=c_k(a,b),&\ \forall (a,b)\in\supp\pi_k\\
    \end{array}\right.,\quad \forall k\geq 1
\end{equation*}
with $\etat_k(a,b)=\langle\eta_k,\epsilon_{(a,b)}\rangle,$
$(a,b)\in\s_k.$ As $\cup_k\s_k=\{c<\infty\},$ this is equivalent
to
\begin{equation*}
    \left\{\begin{array}{ll}
     \etat(a,b)\leq c(a,b),&\ \forall (a,b)\in \{c<\infty\}\\
     \etat(a,b)=c(a,b),&\ \forall (a,b)\in\supp\pi\bigcap\{c<\infty\}\\
    \end{array}\right.
\end{equation*}
Finally, one sees with (\ref{eq-49}) that $\etat=\limn\rho_n$ on
$\s.$ This implies that $\etat$ is measurable on $\s$ and
completes the proof of the lemma.
\endproof

During this proof, we have used the following elementary lemmas.

\begin{lemma}\label{res-12} Let $\ell$ be a \positif\ element of $\Cc'.$
For any $(a,b)\in\mathrm{supp\,}\ell,$  there exists a sequence
$\{h_k\}_{k\geq 1}$ of \positif\  continuous bounded functions on
$\AB$ such that
$\lim_{k\rightarrow\infty}h_k.\ell=\epsilon_{(a,b)}$ in with
respect to $\sigma(\Cc',\Cc).$
\end{lemma}
\proof
 To
see this, consider a \decroissant\  sequence  $\{G_k\}_{k\geq 1}$
of neighbourhoods of $(a,b)$ with $\lim_k G_k=\{(a,b)\}$ and
choose $h_k$ such that $\{h_k>0\}\subset G_k$ and $\langle
h_k,\ell\rangle=1,$ this is possible since $\AB$ is a metric
space.
\endproof

\begin{lemma}\label{res-14a}\
For any $\eta\in\Cc^{\prime *},$ the three following statements
are equivalent:
\begin{enumerate}[(i)]
    \item $\Fd(\eta)=0;$
    \item  $\eta\in \overline{\Gamma};$
    \item $\langle\eta,\ell\rangle\leq \Fs_1(\ell),$ for all $\ell\in\Cc'.$
\end{enumerate}
\end{lemma}

\proof The equivalence $\textrm{(i)}\Leftrightarrow\textrm{(ii)}$
is an immediate consequence of (\ref{eq-72}).
\\
Let us prove: $\textrm{(ii)}\Leftrightarrow\textrm{(iii)}.$ Taking
the closure, it is enough to check that for all $u$ in $\CAB$
\begin{equation}\label{eq-126a}
    u\in \Gamma\Leftrightarrow [\langle u,\ell\rangle\leq
    \Fs_1(\ell), \forall \ell\in\Cc'].
\end{equation}
Young's inequality $\langle u,\ell\rangle\leq \F(u)+\Fs(\ell)$ and
$u\in\Gamma\Leftrightarrow\F(u)=0$ for all $u,\ell$ give the
direct implication. For the converse, choosing
$\ell=\varepsilon_{(a,b)}$ in the right-hand side of
(\ref{eq-126a}), one obtains for all $(a,b)\in \s,$ $u(a,b)\leq
\Fs(\varepsilon_{(a,b)}).$ But, $\Fs(\varepsilon_{(a,b)})=c(a,b)$
by Proposition \ref{res-02}-d. This proves (\ref{eq-126a}) and
completes the proof of the lemma.
\endproof

\subsection{Optimal plan : completing the proofs of Theorem \ref{res-MK2} and
\ref{res-MK3}}\label{sec-04}

We are now in position to complete the proofs of these results.

\proof[Proof of Theorem \ref{res-MK2}]
 \boulette{1}
\emph{Sufficient condition}.\ Let $\pi\in P(\mu,\nu)$ be such that
$\IAB c\,d\pi<\infty.$ Let $\varphi$ and $\psi$ satisfy
(\ref{eq-82}). Because of Lemma \ref{res-19}, one obtains that
$\ls\varphi$ and $\ls\psi$ still satisfy (\ref{eq-82}) as well as
 $\ls\varphi\in L_1(A,\mu)$ and $\ls\psi\in L_1(B,\nu).$
Thanks to Lemma \ref{res-13}, there exists some $\omega\in
\overline{\Upsilon}$ (see (\ref{eq-83})) such that
$\langle\omega,(\mu,\nu)\rangle=\IAB c\,d\pi.$ But, with Lemma
\ref{res-14}-c: $T_2^*\omega\in\overline{\Gamma}.$ Therefore, one
can apply Theorem \ref{res-20}-b which insures that $\pi$ is
optimal.
\par\smallskip\noindent
\emph{Necessary condition}.\ Let $\pi$ be an optimal plan. Because
of Theorem \ref{res-20}-b there exists $\omega\in\Yii$ such that
$\eta:=T_2^*\omega\in\overline{\Gamma}$ and (\ref{eq-17}-b) holds.
With Lemma \ref{res-11c}, one sees that $\etat$ defined by
(\ref{eq-71}) satisfies $\etat\leq c$ on $\s$ and
     $\etat=c$ on $\supp\pi\bigcap\{c<\infty\}.$
By (\ref{eq-120}), for all $(a,b)\in\s$ we have
$\etat(a,b)=\widetilde{\omega}_A(a)+\widetilde{\omega}_B(b)$ where
$\widetilde{\omega}_A(a)=\langle \omega_A,\epsilon_a\rangle$ and
$\widetilde{\omega}_B(b)=\langle \omega_B,\epsilon_b\rangle.$ One
concludes the proof, taking $\varphi=\1_{\sA}\widetilde{\omega}_A$
and $\psi=\1_{\sB}\widetilde{\omega}_B$ where $\sA$ and $\sB$ are
the canonical projections of $\s$ on $A$ and $B.$

\Boulette{2} It appears from Lemmas \ref{res-13} and \ref{res-11c}
that the optimal functions $(\varphi,\psi)$ and the optimal linear
form $\omega\in\Yii$ associated with $\pi$ by the KKT condition
(\ref{eq-17}), see Theorem \ref{TLag2}, are related to each other
by
\begin{equation}\label{eq-84}
    \omega(\epsilon_a,\epsilon_b)=\varphi\oplus\psi(a,b), \textrm{ for $\pi$-a.e.
    }(a,b)\in\AB.
\end{equation}
Therefore, (\ref{eq-17}) and (\ref{eq-26}) express the same KKT
condition. If $\varphi$ and $\psi$ are measurable, then they are
integrable by Lemma \ref{res-19}-1 and they solve $(\overline{K})$
by Theorem \ref{res-20}. This proves statement (a). In the general
situation (b), replacing $(\varphi,\psi)$ by
$(\ls\varphi,\ls\psi),$ one concludes similarly by means of Lemma
\ref{res-19}-2.
\endproof

\proof[Proof of Theorem \ref{res-MK3}] By Theorem \ref{res-MK2}
there exist functions $\varphi_1$ and $\psi_1$ satisfying
(\ref{eq-82}). By Lemma \ref{res-19} there exist functions
$\varphi_2$ and $\psi_2$ such that $\varphi_2\in L_1(A,\mu)$ and
$\psi_2\in L_1(B,\nu).$ Now with Lemma \ref{res-13}, one can
extend $(\varphi_2,\psi_2)$ in the sense of (\ref{eq-84}) into
$\omega\in\Yii$ such that $ \omega(\kappa)\leq \Ls(|\kappa|),$
$\forall \kappa\in \Xi.$ But, this is clearly equivalent to $
|\omega(\kappa)|\leq \Ls(|\kappa|),$ $\forall \kappa\in \Xi.$
Applying Lemma \ref{res-11c} and taking
$\varphi=\1_{\sA}\widetilde{\omega}_A$ and
$\psi=\1_{\sB}\widetilde{\omega}_B$ as in the proof of the
necessary condition of Theorem \ref{res-MK2} leads to the desired
result.
\endproof

\section{The proofs of the results of Section \ref{sec:abstractpb}}
\label{sec:proofs}

The results of Section \ref{sec:abstractpb} are a summing up of
Proposition \ref{P1}, Lemma \ref{C4} , Proposition \ref{P3},
Corollary \ref{res-09}, Lemma \ref{L6b}, Proposition \ref{res-06},
Proposition \ref{res-07} and Proposition \ref{res-08}.

We are going to apply the general results of the Lagrangian
approach to the minimization problem \Po\ which are recalled at
Appendix \ref{sec:convexmin}. We use the notations of Appendix
\ref{sec:convexmin}.

\subsection{Preliminary technical results}
Recall that $|u|_\F=\inf\{\alpha>0 ; \F_{\pm}(u/\alpha)\leq 1\}$
with $\F_{\pm}(u)=\max(\F(u),\F(-u)).$ Its associated dual uniform
norm is
$$
|\ell|_\F^*\eqdef \sup_{u, |u|_\F\leq 1}|\ul|,\quad \ell\in\Li
$$ on
$\Li.$ The topological dual space of $(\Li, |\cdot|_\F^*)$ is
denoted by $\Li'.$ It is the topological bidual space of
$(\Ui,|\cdot|_\F).$
\\
Similarly, recall that $|y|_\La=\inf\{\alpha>0 ;
\La_{\pm}(y/\alpha)\leq 1\}$ with
$\La_{\pm}(y)=\max(\La(y),\La(-y)).$ Its associated dual uniform
norm is
$$
|x|_\La^*\eqdef \sup_{y, |y|_\La\leq 1}|\yx|, \quad x\in\Xi
$$
on $\Xi.$ The topological dual space of $(\Xi, |\cdot|_\La^*)$ is
denoted by $\Xi'.$ It is the topological bidual space of
$(\Yi,|\cdot|_\La).$
\\
The adjoint operator $T_1^\sharp$ which appears at Lemma
\ref{L2}-f below is defined as follows. For all $\omega\in\Xi'$
and all $\ell\in\Li,$ $|\langle
T^\sharp_1\omega,\ell\rangle_{\Li^\aldu,\Li}|=|\langle
\omega,T\ell\rangle_{\Xi',\Xi}|$

\begin{lemma}\label{L2}
Let us assume \HF\ and \HT.
\begin{enumerate}[(a)]
    \item $\dom\Fs\subset\Li$ and $\dom\Ls\subset\Xi$
    \item $T(\dom\Fs)\subset\dom\Ls$ and $T\Li\subset\Xi$
    \item $T$ is $\sLU$-$\sXY$-continuous
    \item $T_2^\ast: \Yii\to\Uii$ is $\sYXii$-$\sULii$-continuous
    \item $T_1: \Li\to\Xi$ is  $\NF^*$-$\NL^*$-continuous
    \item $T_1^\sharp\Xi'\subset \Li'$
    \item $T_1^\ast \Yi\subset\Ui$ and $T_1^\ast: \Yi\to\Ui$ is $\sYXi$-$\sULi$-continuous
    \item $T_1: \Li\to\Xi$ is $\sLUi$-$\sXYi$-continuous
\end{enumerate}
\end{lemma}

\proof \boulette{a} For all $\ell\in\LL$ and $\alpha>0,$ Young's
inequality yields: $\ul= \alpha\langle\ell,u/\alpha\rangle\leq
[\F(u/\alpha)+\Fs(\ell)]\alpha,$ for all $u\in\UU.$ Hence, for any
$\alpha>|u|_\F,$ $\ul\leq [1+\Fs(\ell)]\alpha.$ It follows that
$\ul\leq [1+\Fs(\ell)]|u|_\F.$ Considering $-u$ instead of $u,$
one gets
\begin{equation}
  \label{eq-02}
  |\ul|\leq [1+\Fs(\ell)]|u|_\F, \forall u\in\UU, \ell\in\LL.
\end{equation}
It follows that $\dom\Fs\subset\Li.$ One proves
$\dom\Ls\subset\Xi$ similarly.

\Boulette{b}
 Let us consider
$|\cdot|_{\Fs_{\pm}}$ and $|\cdot|_{\Ls_{\pm}}$ the gauge
functionals of the level sets $\{\Fs_{\pm}\leq 1\}$ and
$\{\Ls_{\pm}\leq 1\}.$ It is easy to show that
\begin{equation}
  \label{eq-11}
  \Ls_{\pm}(x)\leq\Fs_{\pm}(\ell), \mathrm{\ for\ all\ } \ell\in\LL
  \mathrm{\ and\ } x\in\XX \mathrm{\ such\ that\ }T\ell=x
\end{equation}
Therefore, $T(\dom\Fs_{\pm})\subset\dom\Ls_{\pm}.$ On the other
hand, by Proposition \ref{Pgauge} (see the Appendix), the linear
space spanned by $\dom\Fs_{\pm}$ is $\dom |\cdot|_{\Fs_{\pm}}$ and
the linear space spanned by $\dom\Ls_{\pm}$ is $\dom
|\cdot|_{\Ls_{\pm}}.$ But, $\dom |\cdot|_{\Fs_{\pm}}=\dom
|\cdot|_\F^*=\Li$ and $\dom |\cdot|_{\Ls_{\pm}}=\dom
|\cdot|_\La^*=\Xi$ by Proposition \ref{Pgauge} again. Hence,
$T\Li\subset\Xi.$

 \Boulette{c}
 To prove that $T$ is continuous, one has to show that
for any $y\in \YY,$ $\ell\in \LL\mapsto \langle y,T\ell\rangle\in
\mathbb{R}$ is continuous. We get $\ell\mapsto\langle
y,T\ell\rangle=\langle T^*y, \ell\rangle$ which is continuous
since \HTi\ is $T^*y\in \UU.$

\Boulette{d} It is a direct consequence of $T\Li\subset\Xi.$ See
the proof of (c).

\Boulette{e} We know by Proposition \ref{Pgauge} that
$|\cdot|_{\Fs_{\pm}}\sim\NF^*$ and $|\cdot|_{\Ls_{\pm}}\sim\NL^*$
are equivalent norms on $\Li$ and $\Xi$ respectively. For all
$\ell\in\Li,$ $|T\ell|_\La^*\leq 2|T\ell|_{\Ls_{\pm}}
=2\inf\{\alpha>0 ;\Ls_{\pm}(T\ell/\alpha)\leq 1\} \leq
2\inf\{\alpha>0 ; \Fs_{\pm}(\ell/\alpha)\leq 1\}.$ This last
inequality follows from (\ref{eq-11}). Going on, we get
$|T\ell|_\La^*\leq 2|\ell|_{\Fs_{\pm}}\leq 4|\ell|_\F^*,$ which
proves that $T_1$ shares the desired continuity property with
$\|T_1\|\leq 4.$

 \Boulette{f}
 Let us take  $\omega\in\Xi'.$ For all $\ell\in\Li,$
$|\langle T^\sharp_1\omega,\ell\rangle_{\Li^\aldu,\Li}|=|\langle
\omega,T\ell\rangle_{\Xi',\Xi}|$ $\leq \|\omega\|_{\Xi'}
|T\ell|_\La^*\leq \|\omega\|_{\Xi'}\|T_1\| |\ell|_\F^*$ where
$\|T_1\|<\infty,$ thanks to (e). Hence, $T^\sharp _1\omega$ stands
in  $\Li'.$

 \Boulette{g}
Let us take  $y\in\Yi.$ We've just seen that $T^\ast _1y$ stands
in $\Li'.$ Let us show that in addition, it is the strong limit of
a sequence in $\UU.$ Indeed, there exists a sequence $(y_n)$ in
$\YY$ such that $\lim_{n\rightarrow\infty} y_n=y$ in $(\Yi,\NL).$
Hence, for all $\ell\in\Li,$ $|\langle T^\ast _1y_n-T^\ast
_1y,\ell\rangle_{\Li^\aldu,\Li}|=|\langle
y_n-y,T\ell\rangle_{\Yi,\Xi}|$ $\leq \|T_1\| |y_n-y|_\La
|\ell|_\F^*$ and $\sup_{\ell\in\Li, |\ell|_\F^*\leq 1}|\langle
T^\ast _1y_n-T^\ast _1y,\ell\rangle|\leq \|T_1\| |y_n-y|_\La$
tends to 0 as $n$ tends to infinity, where $T_1^\ast y_n=T^\ast
y_n$ belongs to $\UU$ for all $n\geq 1$ by \HTi. Consequently,
$T_1^\ast y\in\Ui.$
\\
The continuity statement now follows from (d).

 \Boulette{h}

By (b), $T_1$ maps $\Li$ into $\Xi$ and because of (g): $T^\ast
\Yi\subset\Ui.$ Hence, for all $y\in\Yi,$ $\ell\mapsto \langle
T_1\ell,y\rangle_{\Xi,\Yi}=\langle\ell,T^\ast y\rangle_{\Li,\Ui}$
is $\sLUi$-continuous.
 This completes the proof of Lemma \ref{L2}.
 \endproof

 Let $\F_0^*,$  $\La_0^*$ and $\La_1^*$ be the convex conjugates of
$\F_0,$  $\La_0$ and $\La_1$ for the dual pairings
$\langle\UU,\LL\rangle,$ $\langle\YY,\XX\rangle$ and
$\langle\Yi,\Xi\rangle.$

\begin{lemma}\label{res-05}
Under the hypotheses \HF\ and \HT, we have
\begin{enumerate}[(a)]
    \item $\F_0=\F_1\leq \F$ on $\UU$\qquad (a')\ $\La_0=\La_1\leq\La$ on $\YY$
    \item $\Fs=\F_0^*$ on $\LL$\qquad\qquad (b')\ $\Ls\leq\La_0^*$ on $\XX$
    \item $\Fs=\F_0^*=\F_1^*$ on $\Li$\quad\ (c')\ $\Ls\leq\La_0^*\leq\La_1^*$ on $\Xi$
\end{enumerate}
\end{lemma}
\proof (a) follows directly from Lemma \ref{L2}-a, (a') from (a)
and (b') from (a').
\\
(b) follows from the general fact that the convex conjugates of a
function and its convex \lsc\ regularization match.
\\
Let us show (c). As $\UU$ is a dense subspace of $\Ui,$ we obtain
that the restriction of $\Fs$ to $\Li$  is also the convex
conjugate of $\F$ (restricted to $\Li$) for the dual pairing
$\langle\Ui,\Li\rangle.$ Now, with the same argument as in (b),
this implies that $\Fs=\F_1^*$ on $\Li.$
\\
(c') follows from (a'), the fact that $\YY$ is a dense subset of
$\Yi,$ the weak continuity of $T_1^*$ which is proved at Lemma
\ref{L2}-g and the lower semicontinuity of $\F_1.$
\endproof

\begin{lemma}\label{L6}
Under the hypothesis \HF,
\begin{enumerate}[(a)]
    \item  $\Fs=\Fs_0$ is $\sLU$-inf-compact and
    \item $\Fs_1$ is $\sLUi$-inf-compact.
\end{enumerate}

\end{lemma}
\proof \boulette{b} We first prove that $\Fs_1$ is
$\sLUi$-inf-compact. Recall that we already obtained at
(\ref{eq-02}) that $
  |\ul|\leq [1+\Fs(\ell)]|u|_\F,
$ for all $u\in\UU$ and $\ell\in\LL.$ By completion, one deduces
that for all $\ell\in\Li$ and $u\in\Ui,$ $|\ul|\leq
[1+\Fs_1(\ell)]|u|_\F$ (recall that $\Fs=\Fs_1$ on $\Li,$ Lemma
\ref{res-05}-c.) Hence, $\Fs_1(\ell)\leq A$ implies that
$|\ell|_\F^*\leq A+1.$ Therefore, the level set $\{\Fs_1\leq A\}$
is relatively $\sLUi$-compact.
\\
By construction, $\Fs_1$ is $\sLUi$-\lsc. Hence, $\{\Fs_1\leq A\}$
is $\sLUi$-closed and $\sLUi$-compact.

 \Boulette{a}
As $\Fs=\Fs_0=\Fs_1$ on $\Li$ (Lemma \ref{res-05}-c),
$\dom\Fs\subset\Li$ (Lemma \ref{L2}-a) and $\UU\subset\Ui,$ it
follows from the $\sLUi$-inf-compactness of $\Fs_1$ that
$\Fs=\Fs_0$ is $\sLU$-inf-compact.
\endproof

\subsection{A first dual equality}
In this  section we only consider the basic spaces $\UU,\LL,\YY$
and $\XX.$ Let us begin applying Appendix \ref{sec:convexmin} with
$\langle P,A\rangle=\langle \UU,\LL\rangle$ and
 $\langle B,Q\rangle=\langle \YY,\XX\rangle$ and the topologies are
 the weak topologies $\sLU,$ $\sUL,$ $\sXY$ and $\sYX.$ The
 function to be minimized is
 $f(\ell)=\Fs(\ell)+\delta_{ C}(T\ell),$ $\ell\in\LL$ where
$\delta_C(x)=\left\{
\begin{array}{ll}
0  &\mathrm{if\ }x\in C\\
+\infty &\mathrm{if\ }x\not\in C
\end{array}
\right.$ denotes the convex indicator of $C.$ The perturbation $F$
of $f$ is Fenchel's one:
\begin{displaymath}
  F_0(\ell,x)=\Fs(\ell)+\delta_{ C}(T\ell+x), \ell\in\LL, x\in\XX.
\end{displaymath}
We assume \HTi: $T^\ast \YY\subset\UU,$ so that the duality
diagram is
\begin{equation}
 \begin{array}{ccc}
\Big\langle\ \UU & , & \LL \ \Big\rangle \\
T^\ast  \Big\uparrow & & \Big\downarrow
 T
\\
\Big\langle\ \YY & , & \XX\ \Big\rangle
\end{array}
\tag{Diagram 0}
\end{equation}
The analogue of $F$ for the dual problem is
\begin{displaymath}
  G_0(y,u)\eqdef \inf_{\ell,x}\{\yx-\ul+F_0(\ell,x)\}
=\inf_{x\in C}\yx-\F_0(T^\ast y+u).
\end{displaymath}
The corresponding value functions are
\begin{eqnarray*}
\varphi_0(x) &=& \inf\{\Fs_0(\ell); \ell\in\LL: T\ell\in
C-x\},\quad x\in\XX
\\
 \gamma_0(u)&=&\sup_{y\in\YY}\{\inf_{x\in C}\yx-\F_0(T^\ast
 y+u)\},\quad  u\in\UU.
\end{eqnarray*}
The primal and dual problems are \Po\ and \Do.

\begin{lemma}\label{res-04}
Assuming \HF\ and \HTi, if $C$ is a $\sXY$-closed convex set,
$F_0$ is jointly closed convex on $\LL\times \XX.$
\end{lemma}

\proof
 As $T$ is linear continuous (Lemma \ref{L2}-c) and  $C$ is closed
convex, $\{(\ell,x); T\ell+x\in C\}$ is closed convex in
$\LL\times \XX.$ As $\Fs$ is closed convex on $\LL,$ its epigraph
is closed convex in $\LL\times \mathbb{R}.$  It follows that $\epi
F_0=(\XX\times\epi \Fs)\cap \{(\ell,x); T\ell+x\in C\}$ is closed
convex, which implies that $F_0$ is convex and \lsc. As it is
nowhere equal to $-\infty$ (since $\inf F_0\geq \inf \Fs>-\infty,$
$F_0$ is also a closed convex function.
\endproof

Therefore, assuming that $ C$ is a $\sXY$-closed convex set, one
can apply the general theory of Appendix \ref{sec:convexmin} since
the perturbation function $F_0$ satisfies the assumptions
(\ref{Fcv}) and (\ref{Fcl}).

\begin{proposition}\label{P1}
Let us assume that \HF\ and \HT\ hold. If $ C$ is convex and
$\sXY$-closed, we have the dual equality
\begin{equation}\label{ed0}
  \inf(P)=\sup(D_0)\in[0,\infty].
\end{equation}
In particular,  for all $x$ in $\XX,$ we have the little dual
equality
        \begin{equation}\label{ped0}
            \inf\{\Fs(\ell); \ell\in\LL, T\ell=x\}=\Ls_0(x)\in
            [0,\infty].
        \end{equation}
\end{proposition}

\proof The identity (\ref{ped0}) is a special case of (\ref{ed0})
with $C=\{x\}.$
\\
To prove (\ref{ed0}), we consider separately the cases where
$\inf(P)<+\infty$ and $\inf(P)=+\infty.$

\par\medskip\noindent\textit{Case where $\inf(P)<+\infty.$}\
Thanks to Theorem \ref{TLag1}-b', it is enough to prove that
$\gamma_0$ is \usc\ at $u=0.$ We are going to prove that
$\gamma_0$ is continuous at $u=0.$ Indeed, for all $u\in\UU,$
\[
-\gamma_0(u)=\inf_y\{\F_0(T^*y+u)-\inf_{x\in C}\yx\}\leq
\F_0(u)\leq \F(u)
\]
where the first inequality is obtained taking $y=0.$ The norm
$\NF$ is designed so that $\F_0$ is  bounded above on a
$\NF$-neighbourhood of zero. By the previous inequality, so is the
convex function $-\gamma_0.$ Therefore, $-\gamma_0$ is
$\NF$-continuous on $\icordom (-\gamma_0)\ni 0.$ As it is convex
and $\Li=(\UU,\NF)',$ it is also $\sigma(\UU,\Li)$-\lsc\  and a
fortiori $\sUL$-\lsc\ , since $\Li\subset\LL.$

\par\medskip\noindent\textit{Case where $\inf(P)=+\infty.$}\ Note
that $\sup(D_0)\geq -\F_0(0)=0>-\infty,$  so that we can apply
Theorem \ref{TLag1}-b. It is enough to prove that
\[
\ls \varphi_0 (0)=+\infty
\]
in the situation where $\varphi_0(0)=\inf(P)=+\infty.$ We have
$\ls \varphi_0(0)=\sup_{V\in \mathcal{N}(0)}\inf\{\Fs_0(\ell);
\ell: T\ell\in C+V\}$ where $\mathcal{N}(0)$ is the set of all the
$\sXY$-open neighbourhoods of $0\in \XX.$ It follows that for all
$V\in \mathcal{N}(0),$ there exists $\ell\in \LL$ such that
$T\ell\in C+V$ and $\Fs_0(\ell)\leq\ls \varphi_0(0).$ This implies
that
\begin{equation}\label{eq-36b}
    T(\{\Fs_0\leq\ls\varphi_0(0)\})\cap (C+V)\not =\emptyset,\quad
    \forall V\in \mathcal{N}(0).
\end{equation}
On the other hand, $\inf(P)=+\infty$ is equivalent to: $ T(\dom
\Fs_0)\cap C=\emptyset. $
\\
Now, we prove ad absurdum that $\ls\varphi_0(0)=+\infty.$ Suppose
that $\ls\varphi_0(0)<+\infty.$ Because of $T(\dom \Fs_0)\cap
C=\emptyset,$ we have a fortiori
\[
T(\{\Fs_0\leq\ls\varphi_0(0)\})\cap C=\emptyset.
\]
As $\Fs_0$ is inf-compact (Lemma \ref{L6}-a) and $T$ is weakly
continuous (Lemma \ref{L2}-c), $T(\{\Fs_0\leq\ls\varphi_0(0)\})$
is a $\sXY$-compact subset of $\XX.$ Clearly, it is also convex.
But $C$ is assumed to be closed and convex, so that by Hahn-Banach
theorem, $C$ and $T(\{\Fs_0\leq\ls\varphi_0(0)\})$ are
\emph{strictly} separated. This contradicts (\ref{eq-36b}),
considering open neighbourhoods $V$ of the origin in
(\ref{eq-36b}) which are open half-spaces. Consequently,
$\ls\varphi_0(0)=+\infty.$ This completes the proof of the
proposition.
\endproof

\subsection{Primal attainment and dual equality}
We are going to consider the following duality diagram, see
Section \ref{sec:problems}:

\begin{equation}
 \begin{array}{ccc}
\Big\langle\ \Ui & , & \Li \ \Big\rangle \\
T_1^\ast  \Big\uparrow & & \Big\downarrow
 T_1
\\
\Big\langle\ \Yi & , & \Xi\ \Big\rangle
\end{array}
 \tag{Diagram 1}
\end{equation}
Note that the inclusions $T_1\Li\subset\Xi$ and  $T_1^\ast
\Yi\subset\Ui$ which are stated in Lemma \ref{L2} are necessary to
validate this diagram.
\\
Let $F_1, G_1$ and $\gamma_1$ be the analogous functions to $F_0,$
$G_0$ and $\gamma_0.$ Denoting $\varphi_1$ the primal value
function, we obtain
\begin{eqnarray*}
F_1(\ell,x) &=& \Fs_1(\ell)+\delta_{\Ci}(T_1\ell+x), \ell\in\Li,\quad x\in\Xi\\
G_1(y,u)    &=& \inf_{x\in\Ci}\yx -\F_1(T_1^\ast y+u),\quad
y\in\Yi,
u\in\Ui\quad \textrm{and}\\
\varphi_1(x) &=& \inf\{\Fs_1(\ell); \ell\in\Li: T_1\ell\in
C_1-x\}, \quad x\in\Xi,
\\
\gamma_1(u)&=& \sup_{y\in\Yi}\{\inf_{x\in\Ci}\yx -\F_1(T_1^\ast
y+u)\},\quad u\in\Ui.
\end{eqnarray*}
It appears that the primal and dual problems are \Pi\ and \Di.

\begin{lemma}
\label{C4} Assuming \HF\ and \HT, the problems \Po\ and \Pi\ are
equivalent: they have the same solutions and
$\inf(P)=\inf(P_1)\in[0,\infty].$
\end{lemma}
\proof It is a direct consequence of $\dom\Fs\subset\Li,$
$T\Li\subset\Xi$ and $\Fs=\Fs_1$ on $\Li,$ see Lemma \ref{L2}-a,b
and Lemma \ref{res-05}-c.
\endproof

\begin{proposition}[Primal attainment and dual equality]\label{P3}
Assume that \HF\ and \HT\ hold.
\begin{enumerate}
    \item[(a)] For all $x$ in $\Xi,$ we have the little dual equality
        \begin{equation}\label{ped}
            \inf\{\Fs(\ell); \ell\in\LL, T\ell=x\}=\Ls_1(x)\in
            [0,\infty].
        \end{equation}
\end{enumerate}
Assume that in addition \HC\ holds.
\begin{enumerate}
    \item[(b)]  We have the dual equalities
\begin{align}
   & \inf(P)=\sup(D_1)\in[0,\infty]\label{ed1}\\
   & \inf(P)=\inf(P_1)=\inf_{x\in\Ci}\Ls_1(x)\in [0,\infty]\label{ed1bis}
\end{align}
    \item[(c)] If in addition $\inf(P)<\infty,$ then \Po\ is attained in $\Li.$
    \item[(d)] Let $\lb\in\Li$ be a solution to $(P),$ then
$\xb\eqdef T\lb$ is a solution to \PX\ and
$\inf(P)=\Fs(\lb)=\Ls_1(\xb).$
\end{enumerate}
\end{proposition}

\proof $\bullet$\quad
 We begin with the proof of (\ref{ed1}). As,
$\inf(P)=\inf(P_1)$ by Lemma \ref{C4}, we have to show that
$\inf(P_1)=\sup(D_1).$ We consider separately the cases where
$\inf(P_1)<+\infty$ and $\inf(P_1)=+\infty.$

\par\medskip\noindent\textit{Case where $\inf(P_1)<+\infty.$}\
Because of \HC, $F_1$ is jointly convex and $F_1(\ell,\cdot)$ is
$\sXYi$-closed convex for all $\ell\in\Li.$ As $T_1^\ast
\Yi\subset\Ui$ (Lemma \ref{L2}), one can apply the approach of
Appendix \ref{sec:convexmin} to the duality Diagram 1. Therefore,
by Theorem \ref{TLag1}-b', the dual equality holds  if $\gamma_1$
is $\sULi$-\usc\ at 0. As in the proof of Proposition \ref{P1}, we
have $-\gamma_1(u)\leq\F_1(u),$ for all $u\in\Ui.$ But $\F_1$ is
the $\sULi$-\lsc\  regularization of $\tilde\F(u)=\left\{
  \begin{array}{ll}
\F(u) &\mathrm{if\ }u\in\UU\\
+\infty &\mathrm{otherwise}
  \end{array}
\right.,$ $u\in\Ui$ and $\F$ is bounded above by $1$ on the ball
$\{u\in\UU ; |u|_\F<1\}.$ As $\Li=(\Ui,\NF)',$ $\F_1$ is also the
$\NF$-regularization of $\tilde\F.$ Therefore, $\F_1$ is bounded
above by $1$ on $\{u\in\Ui ; |u|_\F<1\},$ since $\{u\in\UU ;
|u|_\F<1\}$ is $\NF$-dense in $\{u\in\Ui ; |u|_\F<1\}.$ As
$-\gamma_1 (\leq\F_1)$ is convex and bounded above on a
$\NF$-neighbourhood of $0,$ it is $\NF$-continuous on $\icordom
(-\gamma_1)\ni 0.$ Hence, it is $\sULi$-\lsc\  at $0.$

\par\medskip\noindent\textit{Case where $\inf(P_1)=+\infty.$}\
This proof is a transcription of  the second part of the proof of
Proposition \ref{P1}, replacing $T$ by $T_1,$ $C$ by $C_1,$ all
the subscripts 0 by 1 and using the preliminary results:
 $\Fs_1$ is inf-compact (Lemma \ref{L6}) and $T_1$ is weakly
continuous (Lemma \ref{L2}-h). This completes the proof of
(\ref{ed1}).

\noindent $\bullet$\quad The identity (\ref{ped}) is simply
(\ref{ed1}) with $C_1=\{x\}.$

\noindent $\bullet$\quad Let us prove (c). By Lemma \ref{L2}-h,
$T_1$ is $\sLUi$-$\sXYi$-continuous. Since $\Ci$ is
$\sXYi$-closed, $\{\ell\in\Li ; T\ell\in\Ci\}$ is $\sLUi$-closed.
As $\Fs_1$ is $\sLUi$-inf-compact (Lemma \ref{L6}), it achieves
its infimum on the closed set $\{\ell\in\Li ;T\ell\in\Ci\}$ if
$\inf(P_1)=\inf(P)<\infty.$

\noindent $\bullet$\quad Let us prove (\ref{ed1bis}). The dual
equality (\ref{ed1}) gives us, for all $x_o\in\Ci,$
$\inf(P_1)=\sup_{y\in\Yi}\{\inf_{x\in\Ci}\yx-\La_1(y)\} \leq
\sup_{y\in\Yi}\{\langle x_o,y\rangle-\La_1(y)\}=\Ls_1(x_o).$
Therefore
\begin{equation}
  \label{eq-03}
  \inf(P_1)\leq\inf_{x\in\Ci}\Ls_1(x).
\end{equation}
In particular, equality holds instead of inequality if
$\inf(P_1)=+\infty.$ Suppose now that $\inf(P_1)<\infty.$ From
statement (c), we already know that there exists $\lb\in\Li$ such
that $\xb\eqdef T\lb\in\Ci$ and $\inf(P_1)=\Fs(\lb).$ Clearly
$\inf(P_1)\leq\inf\{\Fs(\ell) ; T\ell=\xb,
\ell\in\Li\}\leq\Fs(\lb).$ Hence, $\inf(P_1)=\inf\{\Fs_1(\ell) ;
T\ell=\xb, \ell\in\Li\}.$ By the little dual equality (\ref{ped})
we have $\inf\{\Fs_1(\ell) ; T\ell=\xb, \ell\in\Li\}=\Ls_1(\xb).$
Finally, we have obtained $\inf(P_1)=\Ls_1(\xb)$ with $\xb\in\Ci.$
Together with (\ref{eq-03}), this leads us to the desired
identity: $\inf(P_1)=\inf_{x\in\Ci}\Ls_1(x).$

\noindent $\bullet$\quad Finally, (d) is a by-product of the proof
of (\ref{ed1bis}).
\endproof

The following result is an improvement of Lemma \ref{res-05}-c'.
\begin{corollary}\label{res-09}
We have $\dom\Ls_1\subset\dom\Ls_0,$ $\dom\Ls_1\subset\Xi$ and
    in restriction to $\Xi,$ $\Ls_0=\Ls_1.$
\end{corollary}
\proof The first part is already proved at Lemma \ref{L2}-a.  The
matching $\Ls_0=\Ls_1$ follows from (\ref{ped0}) and (\ref{ped}).
\endproof

\subsection{Dual attainment}
\label{sec:dualatt} We now consider the following duality diagram
\begin{equation}
 \begin{array}{ccc}
\Big\langle\ \Li & , & \Uii \ \Big\rangle \\
T_1 \Big\downarrow & & \Big\uparrow
 T_2^\ast
\\
\Big\langle\ \Xi & , & \Yii\ \Big\rangle
\end{array}
\tag{Diagram 2}
\end{equation}
where the topologies are the respective weak topologies. The
associated perturbation functions are
\begin{eqnarray*}
  F_2(\ell,x) &=&
  \Fs_1(\ell)+\delta_{C_1}(T\ell+x),\quad
  \ell\in\Li, x\in\Xi\\
G_2(\zeta,\omega) &=& \inf_{x\in C_1}\xo-\Fd(T_2^\ast
\omega+\zeta),\quad \zeta\in\Uii, \omega\in\Yii
\end{eqnarray*}
 As $F_2=F_1,$ the primal problem is \Pi\ and its value function
 is $\varphi_1:$
 \begin{equation}\label{eq-104}
    \varphi_1(x) = \inf_{x'\in C_1-x}\Ls_1(x'),\quad x\in\Xi
\end{equation}
where we used (\ref{ped}). The dual problem is \Dii.

Assume that $\inf(P)<\infty.$ We know by Proposition \ref{P3}-d
that \PX\ admits at least a solution $\xb=T\lb$ where $\lb$ is a
solution to \Pi. Let us consider the following new minimization
problem
\begin{equation}
    \textsl{minimize } \Fs_1(\ell) \quad \textsl{subject to }\quad
T\ell=\xb,\quad \ell\in\Li   \tag{$P_1^{\xb}$}
\end{equation}
Of course $\lb$ is a solution to \Pi\ if and only if it is a
solution to \Pxb\ where $\xb=T\lb.$ Since our aim is to derive a
representation formula for $\lb,$ it is enough to build our
duality schema upon \Pxb\ rather than upon \Pi. The associated
perturbation functions are
\begin{eqnarray*}
  F_2^{\xb}(\ell,x) &=&
  \Fs_1(\ell)+\delta_{\{\xb\}}(T\ell+x),\quad
  \ell\in\Li, x\in\Xi\\
G_2^{\xb}(\zeta,\omega) &=& \langle\xb,\omega\rangle-\Fd(T_2^\ast
\omega+\zeta),\quad \zeta\in\Uii, \omega\in\Yii
\end{eqnarray*}
 As $F_2^{\xb}$ is $F_1$ with $C_1=\{\xb\},$ the primal problem is \Pxb\ and its value
 function is
 \begin{displaymath}
   \varphi_1^{\xb}(x) = \Ls_1(\xb-x),\quad x\in\Xi.
 \end{displaymath}
The dual problem is
\begin{equation}
    \textsl{maximize }\quad \langle\xb,\omega\rangle-\La_2(\omega),
  \quad\omega\in\Yii
  \tag{$D_2^{\xb}$}
\end{equation}

\begin{lemma}\label{L6b}
Under the hypotheses \HF\ and \HT,
 $\Ls_1$ is $\sXYi$-inf-compact.
\end{lemma}
\proof  By (\ref{ped}): $\inf\{\Fs_1(\ell); \ell\in\Li,
T_1\ell=x\}=\Ls_1(x)$ for all $x\in\Xi$ (note that $\Fs=\Fs_1$ on
$\Li$ by Lemma \ref{res-05}-c.) As $T_1$ is continuous (Lemma
\ref{L2}-h) and $\Fs_1$ is inf-compact(Lemma \ref{L6}), it follows
that $\Ls_1$ is also inf-compact.
\endproof

\begin{proposition}[Dual attainment]\label{res-06}
Assume that \HF, \HT\ and \HC\ hold.
\begin{enumerate}[(a)]
 \item Suppose that
 \begin{equation}\label{eq-106}
    C\cap\icordom\Ls_1\not=\emptyset.
\end{equation}
    Then the dual problem \Dii\ is attained in $\Yii.$
 \item Suppose that $ C\cap\dom\Ls_1\not=\emptyset.$ Then, $\inf(P)<\infty$ and we know (see
    Proposition \ref{P3}-d) that \PX\
    admits at least a solution. If in addition, there exists a solution
    $\xb$ to \PX\ such that
\begin{equation}\label{CQbis}
    \xb\in\diffdom\Ls_1,
\end{equation}
then the dual problem \Dxb\ is attained in $\Yii.$
\end{enumerate}
\end{proposition}

\proof \boulette{a}
 As $F_2=F_1,$  one can apply the approach of Appendix
\ref{sec:convexmin} to the duality Diagram 2. Let us denote
$\fss_1$ the $\sXYi$-\lsc\  regularization of $\varphi_1$ and
$\fss_2$ its $\sXYii$-\lsc\  regularization. Since $\Xi$ separates
$\Yi,$ the inclusion $\Yi\subset\Yii$ holds. It follows that
$\fss_1(0)\leq\fss_2(0)\leq\varphi_1(0).$ But we have (\ref{ed1})
which is $\fss_1(0)=\varphi_1(0).$ Therefore, one also obtains
$\fss_2(0)=\varphi_1(0)$ which is the dual equality
\begin{equation}
  \label{ed2}
 \inf(P_1)=\sup(D_2)
\end{equation}
and one  can apply Theorem \ref{TLag1}-c which gives
\begin{equation}\label{eq-114}
     \mathrm{argmax}(D_2)=-\partial\varphi_1(0).
\end{equation}
It remains to show that the value function $\varphi_1$  given at
(\ref{eq-104}) is such that
\begin{equation}\label{eq-105}
    \partial\varphi_1(0)\not=\emptyset.
\end{equation}
As the considered dual pairing  $\langle\Xi,\Yii\rangle$ is the
saturated algebraic pairing, for (\ref{eq-105}) to be satisfied,
by the geometric version of Hahn-Banach theorem, it is enough that
$0\in\icordom\varphi_1.$ But this holds provided that the
constraint qualification (\ref{eq-106}) is satisfied.

\Boulette{b} Let us specialize to the special case where
$C_1=\{\xb\}.$ The dual equality (\ref{ed2}) becomes
\begin{equation}\label{ed3}
 \inf(P_1^{\xb})=\sup(D_2^{\xb})
\end{equation}
and (\ref{eq-105}) becomes
$\partial\varphi_1^{\xb}(0)\not=\emptyset$  which is directly
implied by (\ref{CQbis}).
\endproof

\begin{remark}
The dual equality (\ref{ed3}) is
\begin{equation}\label{eq-112}
    \Ls_1=\Ls_2
\end{equation}
 where these convex conjugates are to be taken respectively with
respect to $\langle\Xi,\Yi\rangle$ and $\langle\Xi,\Yii\rangle.$
Denoting $\overline{\La}_1$ and $\overline{\La}_2$  the convex
$\sYXii$-\lsc\ regularizations of $\La_1$ and $\La_2,$
(\ref{eq-112}) implies the identity
\begin{equation}\label{eq-111}
\overline{\La}_1=\overline{\La}_2.
\end{equation}
Usual results about convex conjugation tell us that
$\Ls_1(\xb)=\sup_{\omega\in\Yii}\{\langle\xb,\omega\rangle-\overline{\La}_1(\omega)\}=\sup(D_2^{\xb})$
and the above supremum is attained at $\ob$ if and only if
$\ob\in\partial_{\Yii}\Ls_1(\xb).$ This is the attainment
statement in Proposition \ref{res-06}-b.
\end{remark}

\subsection{Dual representation of the minimizers}

We keep the framework of Diagram 2 and derive the KKT relations in
this situation. The Lagrangian associated with $F_2=F_1$ and
Diagram 2 is for any $\ell\in\Li, \omega\in\Yii,$
\begin{eqnarray*}
  K_2(\ell, \omega)&\eqdef&\inf_{x\in\Xi}\{\xo+\Fs_1(\ell)+\delta_{\Ci}(T\ell+x)\},\\
&=& \Fs_1(\ell)-\langle T\ell,\omega\rangle + \inf_{x\in\Ci}\xo.
\end{eqnarray*}

\begin{proposition}\label{res-07}
Assume that \HF, \HT\ and \HC\ hold.
\begin{enumerate}[(a)]
 \item Any $\lb\in\Li$ is a solution to \Pi\
 if and only if there exist some $\ob\in\Yii$ such that the following three statements hold
    \begin{enumerate}[(1)]
        \item $T\lb\in C$
        \item $\langle T\lb,\ob\rangle\leq \langle x,\ob\rangle$ for
        all $x\in C_1$
        \item and the following representation formula holds
        \begin{equation}\label{eq-110}
    \lb\in\partial_{\Li}\Fd(T^*_2\ob)
\end{equation}
    \end{enumerate}
    More, these three statements hold if and only if: $\lb$ is
    solution to \Pi, $\ob$ is a solution to \Dii\ and $\inf(P_1)=\sup(D_2).$
    \\
Statement (\ref{eq-110}) is equivalent to the Young's identity
\begin{equation}\label{eq-107}
    \Fs(\lb)+\Fd(T^*_2\ob)=\langle T\lb,\ob\rangle.
\end{equation}
 \item (Assumption \HC\ is useless here).
 Any $\lb\in\Li$ is a solution to $(P^{\xb})$ if and only
     if $T\lb=\xb$ and
     there exists some $\ob\in\Yii$ such that (\ref{eq-110})
     or equivalently (\ref{eq-107}) holds.
     \\
     More, this occurs if and only if: $\lb$ is a
     solution to \Po, $\ob$ is a solution to \Dxb\ with $\xb:=T\lb$ and
    $\inf(P^{\xb})=\sup(D_2^{\xb}).$
\end{enumerate}
\end{proposition}

\proof This proof is an application of Theorem \ref{TLag2}. Under
the general assumptions  \HF, \HT\ and \HC, we have seen at
Proposition \ref{res-06} that the dual equalities (\ref{ed2}) and
(\ref{ed3}) hold true. In both situations (a) and (b), $(\lb,\ob)$
is a saddle-point; all we have to do is to translate the KKT
relations (\ref{KTa}) and (\ref{KTb}).

\Boulette{a} With $K_2$ as above,  (\ref{KTa}) and (\ref{KTb}) are
$\partial_\ell K_2(\lb,\ob)\ni 0$ and
$\partial_\omega(-K_2)(\lb,\ob)\ni 0.$ Since $-\langle
T\ell,\omega\rangle$ is locally weakly upper bounded as a function
of $\omega$ around $\ob$ and as a function of $\ell$ around $\lb,$
one can apply (Rockafellar, \cite{Roc74}, Theorem 20) to derive
$\partial_\ell K_2(\lb,\ob)=\partial\Fs(\lb)-T_2^\ast \ob$ and
$\partial_\omega(-K_2)(\lb,\ob)= \partial(-\inf_{x\in\Ci}\langle
x,\cdot\rangle)+T\lb.$ Therefore the KKT relations are
\begin{eqnarray}
  T_2^\ast \ob&\in &\partial\Fs(\lb)\label{eq-108}\\
 -T\lb &\in & \partial(\delta_{-\Ci}^*)(\ob)\label{eq-109}
\end{eqnarray}
where $\delta_{-\Ci}^*$ is the convex conjugate of the convex
indicator of $-\Ci.$
\\
As a convex conjugate, $\Fs$ is a closed convex functions. Its
convex conjugate is $\Fd.$ Therefore (\ref{eq-108}) is equivalent
to the following equivalent statements
\begin{eqnarray*}
  &&\lb \in \partial\Fd(T_2^\ast \ob)\\
&& \Fs(\lb)+ \Fd(T_2^\ast \ob)=\langle \lb,T_2^*\ob\rangle
\end{eqnarray*}
Similarly, as a convex conjugate $\delta_{-\Ci}^*$ is a closed
convex functions. Its convex conjugate is $\delta_{-\bar\Ci}$
where $\bar\Ci$ stands for the $\sXYii$-closure of $\Ci.$ Of
course, as $\Ci$ is $\sXYi$-closed by hypothesis \HC, it is a
fortiori $\sXYii$-closed, so that $\bar\Ci=\Ci.$ Therefore
(\ref{eq-109}) is equivalent to
\begin{equation}\label{eq-13}
   \delta_{\Ci}(T\lb)+\delta_{-\Ci}^*(\ob)=\langle
-T\lb,\ob\rangle.
\end{equation}
It follows from (\ref{eq-13}) that $\delta_{\Ci}(T\lb)<\infty$
which is equivalent to $ T\lb\in\Ci.$
\\
Now (\ref{eq-13}) is $-\langle
T\lb,\ob\rangle=\delta_{-\Ci}^*(\ob)=-\inf_{x\in\Ci}\langle
x,\ob\rangle$ which is $\langle T\lb,\ob\rangle
=\inf_{x\in\Ci}\langle x,\ob\rangle.$ This completes the proof of
(a).

\Boulette{b} This follows directly from (a) with $\xb=T\lb$ and
$C_1=\{\xb\}.$
\endproof

\begin{remark}
Thanks to Proposition \ref{P3}-d, (\ref{eq-107}) leads us to
\begin{equation}\label{eq-113}
    \Ls_1(\xb)+\La_2(\ob)=\langle\xb,\ob\rangle
\end{equation}
for all $\xb\in\dom\Ls_1$ and all $\ob\in\Yii$ solution to \Dxb.
By Young's inequality:
$\Ls_2(\xb)+\overline{\La}_2(\ob)\geq\langle\xb,\ob\rangle$ and
the identities (\ref{eq-112}), (\ref{eq-113}), we see that
$\overline{\La}_2(\ob)\geq\La_2(\ob).$ But, the reversed
inequality always holds true. Therefore, we have
$
    \overline{\La}_2(\ob)=\La_2(\ob).
$ This proves that
$
    \La_2=\overline{\La}_2 \textrm{ on } \dom\La_2.
$
\end{remark}

\begin{proposition}\label{res-08}
Assume that \HF, \HT\ and \HC\ hold. Any solution $\ob$ of \Dii\
or \Dxb\ shares the following properties
\begin{itemize}
 \item[(a)] $\ob$ stands in the $\sYXx$-closure of $\dom\La_1.$
 \item[(b)] $T_2^\ast \ob$ stands in the $\sULii$-closures of
 $T_1^\ast(\dom\La_1)$ and $\dom\F.$
 \item[(c)] For any $x_o\in\Xi,$ let us denote $j_{D_{x_o}}$ and
    $j_{-D_{x_o}}$ the gauge functionals on $\Xi$ of the convex sets $D_{x_o}$ and
    $-D_{x_o}$ where $D_{x_o}=\{x\in\Xi;\Ls_1(x_o+x)\leq\Ls_1(x_o)+1\}.$
        \begin{itemize}
            \item[-] Let $\ob$ be any solution of \Dii. Then, for any $x_o$ in $
    C\cap\icordom\Ls_1,$ $\ob$ is $j_{D_{x_o}}$-\usc\ and $j_{-D_{x_o}}$-lower semicontinuous at $0.$
            \item[-] Let $\ob$ be any solution of \Dxb\ with $\xb\in\icordom\Ls_1$. Then, $\ob$ is $j_{D_{\xb}}$-\usc\             and $j_{-D_{\xb}}$-lower semicontinuous at $0.$
        \end{itemize}
\end{itemize}
\end{proposition}
\proof
 \boulette{a}
Because of (\ref{eq-113}), we have $\ob\in\dom\La_2.$ As
$\overline{\La}_2\leq \La_2$ and
$\overline{\La}_1=\overline{\La}_2$ (see (\ref{eq-111})), we
obtain $\ob\in\dom\overline{\La}_1$ which implies (a).

 \Boulette{b}
It follows from (a) and the continuity of $T^*_2,$ see Lemma
\ref{L2}-d that $T_2^\ast \ob$ is in the $\sULii$-closure of
 $T_1^\ast(\dom\La_1).$ On the hand, $T_2^\ast \ob\in\dom\Fd$ and
 $\Fd$ is the $\sULii$-closed convex closure of $\F.$ It follows
 that $T_2^\ast \ob$ is in the $\sULii$-closure of $\dom\F.$

 \Boulette{c}
Let $\ob\in$ argmax$(D_2).$ By (\ref{eq-104}) and (\ref{eq-114}),
for all $x\in\Xi$ and any $x_o\in\Ci,$
$\langle-\ob,x\rangle\leq\varphi_1(x)-\varphi_1(0)
\leq\Ls_1(x_o-x)-\varphi(0)\leq\Ls_1(x_o-x).$ It follows that
$\langle \ob,x\rangle\leq\Ls_1(x_o)+1$ for all $x\in D_{x_o}.$
This implies that for all $x\in\Xi,$ $\langle \ob,x\rangle
\leq[1+\Ls_1(x_o)]j_{D_{x_o}}(x).$ Since $j_D(-x)=j_{-D}(x),$ we
finally obtain
\begin{equation*}
  -[1+\Ls_1(x_o)]j_{-D_{x_o}}(x)\leq \langle \ob,x\rangle
\leq [1+\Ls_1(x_o)]j_{D_{x_o}}(x), \forall x\in\Xi
\end{equation*}
for any $x_o\in C_1,$ which is the desired result. Choosing $x_o$
in $\Ci\cap\icordom\Ls_1$ implies that $j_{D_{x_o}}$ is a
nondegerate homogeneous functional.
\\
The second case where $\ob\in$ argmax\Dxb\ is a specialization of
the previous one.
\endproof

\appendix
\section{A short reminder about convex minimization}
\label{sec:convexmin}

To quote easily and precisely some well-known results of convex
minimization while proving our abstract results at Section
\ref{sec:proofs}, we give a short overview of the approach to
convex minimization problems by means of conjugate duality as
developed in Rockafellar's monograph \cite{Roc74}. For complete
proofs of these results, one can also have a look at the author's
lecture notes \cite{Leo05d}.

Let $A$ be a vector space and $f: A\rightarrow [-\infty,+\infty]$
an extended real convex function. We consider the following convex
minimization problem
\begin{equation}
  \textsl{minimize } f(a), a\in A
\tag{$\mathcal{P}$}
\end{equation}
Let $Q$ be another vector space. The perturbation of the objective
function $f$ is a function $F:A\times Q\rightarrow
[-\infty,+\infty]$ such that for $ q=0\in Q,$ $F(\cdot,
0)=f(\cdot).$ The problem $(\mathcal{P})$ is imbedded in a
parametrized family of minimization problems
\begin{equation}
  \textsl{minimize } F(a,q), a\in A
\tag{$\mathcal{P}_q$}
\end{equation}
The value function of $(\mathcal{P}_q)_{q\in Q}$ is
\begin{displaymath}
  \varphi(q)\eqdef \inf(\mathcal{P}_q) =\inf_{a\in
  A}F(a,q)\in[-\infty,+\infty], q\in Q.
\end{displaymath}
Let us assume that the perturbation is chosen such that
\begin{equation}
  \label{Fcv}
  F \mathrm{\ is\ jointly\ convex\ on\ } A\times Q.
\end{equation}
Then, $(\mathcal{P}_q)_{q\in Q}$ is a family of convex
minimization problems and the value function $\varphi$ is convex.

Let $B$ be a vector space in dual pairing with $Q.$ This means
that $B$ and $Q$ are locally convex topological vector spaces in
separating duality such that their topological dual spaces $B'$
and $Q'$ satisfy $B'=Q$ and $Q'=B$ up to some isomorphisms. The
Lagrangian associated with the perturbation $F$ and the duality
$\langle B,Q\rangle$ is
\begin{equation}
\label{KF}
  K(a,b)\eqdef \inf_{q\in Q}\{\bq+F(a,q)\}, a\in A, b\in B.
\end{equation}
 Under (\ref{Fcv}), $K$ is a convex-concave function.
Assuming in addition that $F$ is chosen such that
\begin{equation}
  \label{Fcl}
 q\mapsto F(a,q) \mathrm{\ is\ a\ closed\ convex\ function\
for\ any\ } a\in A,
\end{equation}
one can reverse the conjugate duality relation (\ref{KF}) to
obtain
\begin{equation}
\label{FK}
  F(a,q)=\sup_{b\in B}\{K(a,b)-\bq\}, \forall a\in A, q\in Q
\end{equation}

Introducing another vector space $P$ in separating duality with
$A$ we define the function
\begin{equation}
\label{GK}
  G(b,p)\eqdef \inf_{a\in A}\{K(a,b)-\ap\}, b\in B, p\in P.
\end{equation}
This formula is analogous to (\ref{FK}). Going on symmetrically,
one interprets $G$ as the concave perturbation of the objective
concave function
\begin{displaymath}
  g(b)\eqdef G(b,0), b\in B
\end{displaymath}
associated with the concave maximization problem
\begin{equation}
  \textsl{maximize } g(b), b\in B
\tag{$\mathcal{D}$}
\end{equation}
which is  the dual problem of $(\mathcal{P}).$ It is imbedded in
the family of concave maximization problems $(\mathcal{D}_p)_{p\in
P}$
\begin{equation}
  \textsl{maximize } G(b,p), b\in B
\tag{$\mathcal{D}_p$}
\end{equation}
whose value function is
\begin{displaymath}
  \gamma(p)\eqdef \sup_{b\in B} G(b,p), p\in P.
\end{displaymath}
Since $G$  is jointly concave, $\gamma$ is also concave. We have
the following diagram
\begin{displaymath}
 \begin{array}{lrccclr}
       &           & \gamma(p)  &         & f(a)   &            &\\
       &\Big\langle& P          & ,       & A      &\Big\rangle &\\
G(b,p) &           &            & K(a,b)  &        &            & F(a,q)\\
       &\Big\langle& B          & ,       & Q      &\Big\rangle &\\
       &           & g(b)       &         & \varphi(q)&         &
\end{array}
\end{displaymath}

The concave conjugate of the function $f$ with respect to the dual
pairing $\langle Y,X\rangle$ is $f^\cav(y)=\inf_x\{\langle
y,x\rangle-f(x)\}$ and its superdifferential at $x$ is
$\widehat{\partial}f(x)=\{y\in Y; f(x')\leq f(x)+\langle y,
x'-x\rangle\}.$

\begin{theorem}
\label{TLag1} We assume that $\langle P,A\rangle$ and $\langle
B,Q\rangle$ are topological dual pairings.
\begin{itemize}
 \item[(a)] We have $\sup(\mathcal{D})=\varphi^{**}(0).$ Hence, the dual equality
$\inf(\mathcal{P})=\sup(\mathcal{D})$ holds if and only if
$\varphi(0)=\varphi^{**}(0).$
 \item[(b)] In
particular,
\[\left.
\begin{array}{l}
  \bullet\ F \textrm{ is jointly convex }  \\
  \bullet\ \varphi\textrm{ is lower semicontinuous at } 0 \\
  \bullet\ \sup(\mathcal{D}) >-\infty\\
\end{array}\right\}\Rightarrow \inf(\mathcal{P})=\sup(\mathcal{D})
\]
 \item[(c)] If the dual equality
holds, then
\begin{displaymath}
  \mathrm{argmax\ }g=-\partial \varphi(0).
\end{displaymath}
\end{itemize}

Let us assume in addition that (\ref{Fcv}) and (\ref{Fcl}) are
satisfied.
\begin{itemize}
 \item[(a')] We have $\inf(\mathcal{P})=\gamma^{\cav\cav}(0).$
Hence, the dual equality $\inf(\mathcal{P})=\sup(\mathcal{D})$
holds if and only if $\gamma(0)=\gamma^{\cav\cav}(0).$
 \item[(b')] In
particular,
\[\left.
\begin{array}{l}
  \bullet\ \gamma\textrm{ is \usc\ at } 0 \\
  \bullet\ \inf(\mathcal{P})<+\infty\\
\end{array}\right\}\Rightarrow \inf(\mathcal{P})=\sup(\mathcal{D})
\]
 \item[(c')] If the dual equality
holds, then
\begin{displaymath}
  \mathrm{argmin\ }f=-\widehat{\partial}\gamma(0).
\end{displaymath}
\end{itemize}
\end{theorem}

\begin{definition}[Saddle-point]
One says that $(\ab,\bb)\in A\times B$ is a \emph{saddle-point} of
the function $K$ if
\[
K(\ab,b)\leq K(\ab,\bb)\leq K(a,\bb),\quad\forall a\in A, b\in B.
\]
\end{definition}

\begin{theorem}[Saddle-point theorem and KKT relations]\label{TLag2}
The following statements are equivalent.

\begin{enumerate}[(1)]
    \item The point $(\ab,\bb)$ is a saddle-point of the Lagrangian $K$
    \item $f(\ab)\leq g(\bb)$
    \item The following  three statements hold
         \begin{enumerate}[(a)]
         \item we have the \emph{dual equality}: $\sup(\mathcal{D})= \inf(\mathcal{P}),$
         \item $\ab$ is a solution to the primal problem $(\mathcal{P})$ and
         \item $\bb$ is a solution to the dual problem $(\mathcal{D})$.
         \end{enumerate}
\end{enumerate}
In this situation, one also gets
\begin{equation}\label{eq-302}
    \sup(\mathcal{D})=
\inf(\mathcal{P})=K(\ab,\bb)=f(\ab)=g(\bb).
\end{equation}
Moreover, $(\ab,\bb)$ is a saddle-point of $K$ if and only if it
satisfies
\begin{eqnarray}
  \partial_a K(\ab,\bb)&\ni& 0 \label{KTa} \\
 \widehat{\partial}_b  K(\ab,\bb) &\ni& 0  \label{KTb}
\end{eqnarray}
where the subscript $a$ or $b$ indicates the unfixed variable.
\end{theorem}

\section{Gauge functionals associated with a convex function}
\label{sec:gauge}

The following result is well-known, but since I didn't find a
reference for it, I give its short proof.

Let $\theta: S\rightarrow[0,\infty]$ be an extended \positif\
convex function on a vector space $S,$ such that $\theta(0)=0.$
Let $\Ss$ be the algebraic dual space of $S$ and $\ts$ the convex
conjugate of $\theta:$
 \begin{displaymath}
   \ts(r)\eqdef \sup_{s\in S}\{\langle r,s \rangle-\theta(s)\}, r\in\Ss.
 \end{displaymath}
It is easy to show that $\ts: \Ss\rightarrow[0,\infty]$ and
$\ts(0)=0.$ We denote $C_\theta\eqdef\{\theta\leq 1\}$ and
$C_{\theta^*}\eqdef\{\theta^*\leq 1\}$ the unit level sets of
$\theta$ and $\theta^*.$ The gauge functionals to be considered
are
\begin{eqnarray*}
  \jt(s)&\eqdef & \inf\{\alpha>0 ;  s\in \alpha C_\theta\}
=\inf\{\alpha>0 ; \theta(s/\alpha)\leq 1\}\in [0,\infty], s\in S.\\
\jts(r)&\eqdef &\inf\{\alpha>0 ;  r\in\alpha C_{\theta^*}\}
=\inf\{\alpha>0 ; \theta^*(r/\alpha)\leq 1\}\in [0,\infty], r\in
\Ss.
\end{eqnarray*}
As 0 belongs to $\Ct$ and $\Cts,$ one easily proves that $\jt$ and
$\jts$ are positively homogeneous. Similarly, as $\Ct$ and $\Cts$
are convex sets, $\jt$ and $\jts$ are convex functions.
\begin{proposition}
\label{Pgauge} Let $\theta: S\rightarrow[0,\infty]$ be an extended
\positif\  convex function on a vector space $S,$ such that
$\theta(0)=0$ as above. Then for all $r\in\Ss,$ we have
\begin{displaymath}
  \frac{1}{2}\jts(r)\leq \ds(r)\eqdef \sup_{s\in\Ct}\langle r,s\rangle
\leq 2 \jts(r).
\end{displaymath}
We also have
\begin{displaymath}
  \mathrm{cone\ }\dom\ts=\dom\jts=\dom\ds
\end{displaymath}
where $\mathrm{cone\ }\dom\ts$ is the convex cone (with vertex
$0$) generated by $\dom\ts.$
\end{proposition}
\proof $\bullet$\ Let us first show that $\ds(r)\leq 2\jts(r)$ for
all $r\in\Ss.$ If $\jts(r)>0,$ then for all $s\in\Ct,$ $\langle
r,s \rangle =\langle r/\jts(r),s\rangle \jts(r)\leq
[\theta(s)+\ts(r/\jts(r))]\jts(r)\leq (1+1)\jts(r).$
\\
If $\jts(r)=0,$ then $\ts(tr)\leq 1$ for all $t>0.$ For any
$s\in\Ct,$ we get $\langle r,s\rangle =\frac{1}{t}\langle
tr,s\rangle\leq \frac{1}{t}[\theta(s)+\ts(tr)]\leq 2/t.$ Letting
$t$ tend to infinity, one obtains that $\langle r,s\rangle \leq
0.$
\\
$\bullet$\ Let us show that $\jts(r)\leq 2\ds(r).$ If
$\ds(r)=\infty,$ there is nothing to prove. So, let us suppose
that $\ds(r)<\infty.$ As $0\in\Ct,$ we have $\ds(r)\geq 0.$
\\
\emph{First case: } $\ds(r)>0.$\ For all $s\in S$ and
$\epsilon>0,$ we have $s/[\jt(s)+\epsilon]\in\Ct.$ It follows that
$\langle r/\ds(r),s\rangle =\langle r,s/[\jt(s)+\epsilon]\rangle
\frac{\jt(s)+\epsilon}{\ds(r)}\leq \ds(r)
\frac{\jt(s)+\epsilon}{\ds(r)}=\jt(s)+\epsilon.$ Therefore,
$\langle r/\ds(r),s\rangle\leq \jt(s),$ for all $s\in S.$
\\
If $s$ doesn't belong to $\Ct,$ then $\jt(s)\leq\theta(s).$ This
follows from the the assumptions on $\theta:$ convex function such
that $\theta(0)=0=\min\theta$ and the positive homogeneity of
$\jt.$ Otherwise, if $s$ belongs to $\Ct,$ we have $\jt(s)\leq 1.$
Hence, $
  \langle r/\ds(r),s\rangle \leq\max(1,\theta(s)), \forall s\in S.
$ On the other hand, there exists $s_o\in S$ such that
$\ts(r/[2\ds(r)])\leq \langle r/[2\ds(r)],s_o\rangle
-\theta(s_o)+1/2.$ The last two inequalities provide us with
$\ts(r/[2\ds(r)])\leq
\frac{1}{2}\max(1,\theta(s_o))-\theta(s_o)+\frac{1}{2}\leq 1$
since $\theta(s_o)\geq 0.$ We have proved that $\jts(r)\leq
2\ds(r).$
\\
\emph{Second case: } $\ds(r)=0.$\ We have $\langle r,s\rangle\leq
0$ for all $s\in\Ct.$ As $\dom\theta$ is a subset of the cone
generated by $\Ct,$ we also have for all $t>0$ and
$s\in\dom\theta,$ $\langle tr,s\rangle\leq 0.$ Hence $\langle
tr,s\rangle -\theta(s)\leq 0$ for all $s\in S$ and $\ts(tr)\leq
0,$ for all $t\geq 0.$ As $\ts\geq 0,$ we have $\ts(tr)=0,$ for
all $t\geq 0.$ It follows that $\jts(r)=0.$ This completes the
proof of the equivalence of $\jts$ and $\ds.$
\\
$\bullet$\ Finally, this equivalence implies that
$\dom\jts=\dom\ds$ and as $\ts(0)=0$ we have $0\in\dom\ts$ which
implies that $\mathrm{cone\ }\dom\ts=\dom\jts.$
\endproof

\bibliographystyle{plain}

\end{document}